\documentclass[a4paper,10pt,reqno]{amsart}
\usepackage{latexsym,amssymb,amsfonts,mathrsfs,verbatim,fancyhdr,multirow,hyperref}
\theoremstyle{plain}
\newtheorem{theorem}{Theorem}[section]

\newtheorem{lemma}[theorem]{Lemma}
\newtheorem{proposition}[theorem]{Proposition}
\newtheorem{definition}[theorem]{Definition}
\newtheorem{corollary}[theorem]{Corollary}
\newtheorem{remark}[theorem]{Remark}

\DeclareMathOperator{\Hg}{Hg}
\DeclareMathOperator{\Id}{Id}
\DeclareMathOperator{\projection}{Pr}
\DeclareMathOperator{\Ad}{Ad}
\DeclareMathOperator{\Imgn}{Im}

\DeclareMathOperator{\End}{End}
\DeclareMathOperator{\tr}{tr}
\DeclareMathOperator{\Tr}{Tr}

\DeclareMathOperator{\Identity}{Id}

\DeclareMathOperator{\Sym}{Sym}
\DeclareMathOperator{\Graph}{Graph}

\DeclareMathOperator{\Cent}{Center}

\begin{document}
\author{Meng Chen}
\title{Complex multiplication, rationality\\
and mirror symmetry for abelian varieties
}
\address{Mathematisches Institut, 
Universit\"at Bonn,
Beringstrasse 4,
53115, Bonn, Germany}
\email{mchen@math.uni-bonn.de}
\date{\today}
\begin{abstract}
We show that
complex multiplication on abelian varieties 
is equivalent to the existence of a constant rational K\"ahler metric.
We give a sufficient condition
for a mirror of an abelian variety of CM-type to be of CM-type.
We also study the relationship between complex multiplication
and rationality of a toroidal lattice vertex algebra.
\end{abstract}
\maketitle

%\tableofcontents

\section{Introduction and results}
This article is inspired by Gukov-Vafa's paper
\cite{GukovVafaCM},
where they shared their insight on the interplay
between rational conformal field theories (CFTs),
mirror symmetry
and complex multiplication
on Calabi-Yau varieties.
Here we study the case of complex abelian varieties
of arbitrary dimension, 
where the above notions find a sound mathematical foundation.

Complex multiplication (CM) on abelian varieties 
has been extensively studied in geometry as well as in number
theory (see e.g. \cite{ShimuraCM},\cite{LangCM}).
Roughly speaking, abelian varieties of CM-type
have by definition the biggest endomorphism algebra 
(see Definition \ref{Def:CM}).
It is a property which solely depends on the complex structure 
of the abelian variety.
For example, in dimension one,
an elliptic curve 
$E\cong\mathbb{C}/\mathbb{Z}\oplus\tau \mathbb{Z}$
is of CM-type
if and only if $\tau$ lies in an imaginary quadratic number field
$\mathbb{Q}(\sqrt{-D})$.
In physics, toroidal CFTs are also very familiar objects
(see e.g. \cite{WlThesis} and \cite{EnckevortSCFT}
and the references therein).
From a physicist's point of view,
a CFT is called rational,
if its partition function
can be written as a finite sum
of the product of holomorphic and anti-holomorphic characters
(see \cite[\S2]{GukovVafaCM}),
and this is the case when its chiral part is maximal.

The interplay between CM
and rational CFT having a real 2-torus as target space is already known in
Moore's paper \cite[\S10]{MooreArithAttractor}, 
and it is generalized by Wendland
to real tori of arbitrary dimension.
We can rephrase her Theorem 4.5.5
in \cite{WlThesis} as follows:
(a) A real torus $\mathbb{T}$ admits a rational constant metric $G$
if and only if $\mathbb{T}$ can be endowed with a complex structure 
$I$ such that the complex torus $(\mathbb{T},I)$ 
is isogenous to a product of elliptic curves of CM-type.
(b) A CFT $\mathcal C(\mathbb{T},G,B)$
associated to a real torus $\mathbb{T}$
endowed with a constant metric $G$ and 
a B-field $B\in H^2(\mathbb{T},\mathbb{R})$
is rational if and only if both $G$ and $B$ are rational.
(c) Combining (a) and (b) one can say that
a real torus $\mathbb{T}$ admits a rational CFT $\mathcal C(\mathbb{T},G,B)$
if and only if 
$\mathbb{T}$ can be endowed with a complex structure 
$I$ such that $(\mathbb{T},I)$ 
is isogenous to a product of elliptic curves of CM-type.

In \cite{GukovVafaCM} Gukov and Vafa relate this to mirror symmetry.
They make the following observation:
Let $(E'\cong\mathbb{C}/\mathbb{Z}\oplus\rho \mathbb{Z},G',B')$
be a mirror elliptic curve of 
$(E\cong\mathbb{C}/\mathbb{Z}\oplus\tau \mathbb{Z},G,B)$,
where $G$ is a constant K\"ahler metric.
Then the N=2 CFT $\mathcal C(E,G,B)$
is rational if and only if both $E'$ and $E$ are of CM-type
over the same imaginary quadratic field, 
i.e. $\tau,\rho\in \mathbb{Q}(\sqrt{-D})$,
in this case, $E'$ and $E$ are in particular isogenous.
Note that the difference between this statement and
Wendland's result is that Gukov and Vafa start with an 
elliptic curve with a \textit{given} complex structure
endowed with a \textit{K\"ahler} metric,
and use an additional N=2 structure on the CFT to encode the complex structure
of the target space, and hence link the complex geometry (e.g. CM)
of the mirror pair.

One then naturally asks whether
similar relations hold for abelian varieties
of arbitrary dimension,
where, along with CM and mirror symmetry,
also CFTs and non-linear sigma model received a mathematical treatment,
this is given by Kapustin and Orlov
in terms of vertex algebras (see \cite{KO}).
Note that in their work, it is already apparent that 
mirror symmetry of the target space, i.e. the complex tori,
is equivalent to mirror symmetry of the N=2 vertex algebras
(see Theorem 4.5 in \cite{KO}).
We shall use a slightly different construction, which we call
a \textit{lattice vertex algebra} (for in a special case, it reduces
to the the well-known lattice vertex algebra constructed in
\cite[\S5.4]{KacVA}).
It turns out that it is isomorphic to the Kapustin-Orlov vertex algebra
(see Appendix),
but more suitable for our discussion of CM and rationality.
Given these properly defined notions, we formulate the question
which we fully answer in this article:

\smallskip
\begin{center}
  \begin{tabular}{p{10cm}lp{3cm}r}
\textit{Let $(X',G',B')$ be a mirror partner of $(X,G,B)$,
where $X$ is an abelian variety endowed with
a constant K\"ahler metric $G$
and a B-field in $H^2(X,\mathbb{R})$.
Is the N=2 lattice vertex algebra $V(X,G,B)$ rational
if and only if $X$ and $X'$ are isogenous and both of CM-type?}
&\multirow{4}{*}{\qquad(Q)}
  \end{tabular}
\end{center}
\smallskip

Let us now explain our main results.
The following theorem shows that CM, which is \textit{a priori}
determined solely
by the complex structure, turns out to be equivalent to the rationality of 
a K\"ahler metric.

\smallskip\noindent
\textbf{Theorem \ref{Prop:RatKMetIffCM}.}
\textit{An abelian variety $X$ is of CM-type
if and only if $X$ admits a constant rational K\"ahler metric.}

\smallskip\noindent
(A rational K\"ahler metric is
a K\"ahler metric which takes only rational values on the lattice
$\Gamma$ of $X$.)
This theorem is apparently independent of mirror symmetry and CFT,
but it plays an important r\^ole in the 
relation between CM and the rationality of N=2 lattice vertex algebras.
This will be evident through the following theorem.

\smallskip\noindent
\textbf{Theorem \ref{Thm:VARatGBRat}.}
\textit{The N=2 lattice vertex algebra 
$V(T,G,B)$
associated to a complex torus $T$ endowed with
a constant K\"ahler metric $G$ and a B-field
is rational if and only if $G$ and $B$ 
are both rational.}

\smallskip\noindent
A few comments are due here.
The rationality of a N=2 lattice vertex algebra 
is defined on the underlying lattice vertex algebra 
(i.e. without the N=2 structure, see Definition \ref{Def:Rational}).
Hence forgetting the N=2 structure and the complex structure of $T$,
this theorem is in complete accordance with Wendland's result
mentioned earlier.
Combining the last two theorems we have

\smallskip\noindent
\textbf{Corollary \ref{Prop:CMVAlgRat}.} 
\textit{
An abelian variety $X$ is of CM-type
if and only if 
$X$ admits a rational N=2 lattice vertex algebra
$V(X,G,B)$.
}

\smallskip
Now we make the link between CM and mirror symmetry for tori.
In Section \ref{Sect:CMMirror}
we define mirror symmetry in terms of generalized K\"ahler structures
(GKS) (see Definition \ref{Def:gkstr})
which have the advantage of treating the complex 
and the symplectic structures of a triple $(T,G,B)$ on equal footing.
For $(T,G,B)$ one can define the following GKS 
$(\mathcal{I},\mathcal{J})$:
\begin{equation*}
  \mathcal{I}:
=\left(
    \begin{matrix}
      I&0\\
      BI+I^{t}B&-I^t
    \end{matrix}
\right)
\quad\mbox{and}\quad
  \mathcal{J}:
=\left(
    \begin{matrix}
      \omega^{-1}B&-\omega^{-1}\\
      \omega+B\omega^{-1}B&-B\omega^{-1}
    \end{matrix}
\right),
\end{equation*}
where $\omega$ is the K\"ahler form $G(\cdot,I\cdot)$
(see Definition \ref{Def:IJinduced} for more details).
Note that the composition $\mathcal{I}\mathcal{J}$
of such an induced GKS is defined over $\mathbb{Q}$ if and only if 
both $G$ and $B$ are rational (see 
Definition \ref{Def:ratIJ} and Lemma \ref{Lem:LatticeDecomposition} (iv)).
A pair $(T,\mathcal{I},\mathcal{J})$ and $(T',\mathcal{I}',\mathcal{J}')$
is called a mirror pair if there is a mirror map which
exchanges $\mathcal{I}$ with $\mathcal{J}'$ and $\mathcal{J}$
with $\mathcal{I}'$ (see Definition \ref{Def:MS}).
We find a sufficient condition for CM to be transmitted to mirror partners:

\smallskip\noindent
\textbf{Theorem \ref{Prop:CMMirror}.}
\textit{Let $(X,G,B)$
and $(X',G',B')$
be mirror abelian varieties.
Suppose $X$ is of CM-type.
If both $G$ and $B$ are rational,
then $X$ and $X'$ are isogenous.
In particular, $X'$ is also of CM-type.}

\smallskip\noindent
The converse of this theorem however does not hold:

\smallskip\noindent
\textbf{Proposition \ref{Prop:MirrorNotRat}.}
\textit{There are mirror abelian varieties
$(X,G,B)$ and $(X',G',B')$,
such that $X$ and $X'$ are isogenous and of CM-type,
but neither $\mathcal{I} \mathcal{J}$
nor $\mathcal{I}' \mathcal{J}'$
is defined over $\mathbb{Q}$,
where $(\mathcal{I},\mathcal{J})$
and $(\mathcal{I}',\mathcal{J}')$
denote their induced GKS.}

\smallskip

Combining all the results from the above,
we now formulate our answer to the question (Q):

\smallskip\noindent
\textbf{Corollary \ref{Prop:Summary}.} 
\textit{Let $(X,G,B)$
and $(X',G',B')$ be
mirror abelian varieties.
If the N=2 lattice vertex algebra $V(X,G,B)$ 
is rational,
then $X$ and $X'$ are isogenous and both of CM-type.
Conversely, \,however, 
\,there exist mirror abelian varieties
$(X,G,B)$ and
$(X',G',B') $
such that $X$ and $X'$ are isogenous and both of CM-type,
but neither $V(X,G,B)$ 
nor $V(X',G',B')$ is rational.
}

\smallskip

This paper is organized as follows:
Section \ref{Sect:CMRatMet} presents the proof of 
Theorem \ref{Prop:RatKMetIffCM} which is a pure geometric result.
Section \ref{Sect:CMMirror} formulates mirror symmetry for tori
and we prove Theorem \ref{Prop:CMMirror}.
In Section \ref{Sect:CounterExp} we give an explicit 
counter-example to the converse of Theorem \ref{Prop:CMMirror},
proving Proposition \ref{Prop:MirrorNotRat}.
Section \ref{Sect:VA} deals with vertex algebras.
We construct the lattice vertex algebra,
which is isomorphic to the Kapustin-Orlov vertex algebra (shown in Appendix),
and we define the notion of rationality,
and we prove Theorem \ref{Thm:VARatGBRat}, 
Corollary \ref{Prop:CMVAlgRat} and 
Corollary \ref{Prop:Summary}.

\medskip\noindent
\textbf{Acknowledgments:} 
This paper is a part of my Ph.-D. thesis.
I am very grateful to my thesis advisor Prof. D. Huybrechts
for his great help.
It is also a pleasure to thank
Prof. B. van Geemen, Prof. M. Gaberdiel 
and Mark Rosellen 
for corrections.
Discussions with and comments by Christian van Enckevort,
Katrin Wendland, Oren Ben-Bassat and Gregory Moore
were also very instructive.

%\bigskip
\section{Complex multiplication and rational K\"ahler metric}\label{Sect:CMRatMet}
The aim of this section is to prove Theorem \ref{Prop:RatKMetIffCM}.
Let us first explain when a constant K\"ahler metric is called 
rational.
If we identify the tangent space of a complex torus 
$T=\mathbb{C}^g/\Gamma$
with $\Gamma_{\mathbb{R}}:=\Gamma\otimes_{\mathbb{Z}} \mathbb{R}$,
then the complex structure $I$ of $T$ can be considered as an endomorphism
of $\Gamma_{\mathbb{R}}$ with $I^2=-\Identity$.
A constant K\"ahler metric $G$ is a positive definite bilinear form on 
$\Gamma_{\mathbb{R}}\times\Gamma_{\mathbb{R}}$,
which is compatible with $I$,
i.e. $G(I\cdot,I\cdot)=G(\cdot,\cdot)$.
If $G$ takes only rational values on $\Gamma\times\Gamma$,
then we call it \textit{rational}.
For complex multiplication we adopt Mumford's definition
(\cite[\S2]{MumfordShimura}):
\begin{definition}\label{Def:CM}
A simple abelian variety $X$ of dimension $g$ is \textup{of CM-type over} $K$
if there is an embedding 
$K\hookrightarrow\End_{\mathbb{Q}}(X):=\End(X)\otimes_{\mathbb{Z}}\mathbb{Q}$
of an algebraic number field $K$ of degree $2g$ over $\mathbb{Q}$
into the endomorphism algebra of $X$.
More generally, an (not necessarily simple) abelian variety $X$
is \textup{of CM-type} if $X$ is isogenous to a product
of simple abelian varieties of CM-type.
\end{definition}
One can show (see \cite[\S5]{ShimuraCM},\cite[\S2]{MumfordShimura})
\begin{proposition}
  (i) If a \textup{simple} abelian variety $X$ of dimension $g$
is of CM-type over $K$,
then $K$ is necessarily a CM-field
(i.e. a totally complex quadratic extension of 
a totally real number field)
of rank $2g$
over $\mathbb{Q}$.

(ii) More generally, an (not necessarily simple) 
abelian variety $X$ of dimension $g$ is of CM-type if and only if 
$\End_{\mathbb{Q}}(X)$ contains a commutative semi-simple algebra
of rank $2g$ over $\mathbb{Q}$.

(iii) If $\End_{\mathbb{Q}}(X)$ of an abelian variety $X$ of dimension $g$
contains a number field of degree $2g$ over $\mathbb{Q}$,
then $X$ is isogenous to a product $B\times\cdots\times B$
with a simple abelian variety $B$ of CM-type.
\end{proposition}

\begin{remark}
\textup{Definition \ref{Def:CM} 
is stronger than the notion of admitting complex multiplication 
in \cite{GukovVafaCM}.
There a complex torus 
$T=\mathbb{C}^g/(\mathbf{1}\,\,\mathcal T)\mathbb{Z}$
is said to ``admit complex multiplication''
if there is a non-trivial endomorphism $A\in GL(g,\mathbb{C})$,
such that there exist integer matrices $M',N',M,N$ and
$N$ has rank $g$ and
$
A=M+N \mathcal{T} 
\quad\mbox{and}\quad
\mathcal{T} A=M'+N' \mathcal{T}$.
For example, if $\End_{\mathbb{Q}}(X)$ contains a number field 
$\mathbb{Q}(\xi)$ 
of rank $2g=2\dim X$,
then an integral multiple of the multiplication by $\xi$
on $\mathbb{C}^g$ would satisfy Gukov-Vafa's condition.
However, one can show that the converse does not hold.}
\end{remark}

Next we recall a few general facts
which we will need in the sequel,
e.g. the construction of simple abelian varieties
of CM-type and the characterization of complex multiplication
by the Hodge group.

It can be shown 
that any simple abelian variety of CM-type over a CM-field  $K$ is isogenous
to an abelian variety constructed in the following way
(see \cite[\S6]{ShimuraCM} or \cite[Chap. IV]{MumAV}).
Let $K$ be a quadratic extension of 
a real field $K_0$
and denote by $\mathcal O_K$ its ring of integers.
Let $\Phi=\{\sigma_1,\ldots,\sigma_g\}$
be a CM-type,
i.e. $\Phi\cup\bar{\Phi}$ is the whole set of the embeddings of $K$ into $\mathbb{C}$.
One obtains a complex torus
$\mathbb{C}^g/\Phi(\mathcal O_K)$,
where $\Phi(\mathcal O_K)$ is the image of the embedding of $\mathcal O_K$
into $\mathbb{C}^g$ given by $x\mapsto(\sigma_1x,\ldots,\sigma_gx)$.
There is always a Riemann form, which makes this complex torus
into an abelian variety.
Indeed, there exists an element $\beta\in\mathcal O_K$
such that $K=K_0(\beta)$
and $-\beta^2$ is totally positive and
$\Imgn\sigma_j(\beta)>0,\forall\sigma_j\in \Phi$.
The Riemann form on the tangent space is defined as
\begin{equation}\label{Eq:RiemForm}
  E(z,w)=\sum_{j=1}^g\sigma_j(\beta)\left(
\sigma_j(z)\overline{\sigma_j(w)}-\overline{\sigma_j(z)}\sigma_j(w)
\right).
\end{equation}
On $\mathcal O_K$, it takes values in $\mathbb{Q}$:
\[
E(a_1,a_2)
=\Tr_{K/\mathbb{Q}}(\beta a_1\bar a_2)
\in \mathbb{Q},
\qquad \forall a_1,a_2\in \mathcal O_K.
\]
The thus-constructed abelian variety 
is simple if and only if there is \textit{no}
proper subfield $L$ of $K$
with the properties
(a) $L$ is a purely complex quadratic 
extension of $L\cap K_0$
and (b) if $\sigma_i|_{L\cap K_0}=\sigma_j|_{L\cap K_0}$
then $\sigma_i|_{L}=\sigma_j|_{L},
\forall \sigma_i,\sigma_j\in\Phi$.

As to the Hodge group,
it is one of the main tools to study Hodge structures.
In the case of abelian varieties,
the Hodge structure is of weight one
and the Hodge group is defined to be the smallest algebraic
subgroup of $GL(\Gamma_{\mathbb{Q}})$
defined over $\mathbb{Q}$, 
whose $\mathbb{R}$-points include the image of the unit circle $S^1$ 
under the map
\[
h: S^1\longrightarrow SL(\Gamma_{\mathbb R}),
\qquad
x+yi\longmapsto x+yI
\]
i.e. $h(S^1)\subset\Hg(X)(\mathbb{R})$.
Moreover, Hodge groups are reductive and connected (see \cite{MumfordFamAV}).
There is the following characterization of the complex multiplication.
The equivalence between
(ii) and (iii) below is also well known to the experts.
Nevertheless, since we could not find an explicit proof in the literature,
we include a complete proof for the readers' convenience.
\begin{proposition}\label{Prop:HgCpcComm}
Let $X$ be an abelian variety.
Then the following conditions are equivalent:
\begin{itemize}
\item[(i)] $X$ is of CM-type.
\item[(ii)] $\Hg(X)$ is commutative.
\item[(iii)] $\Hg(X)(\mathbb{R})$ is compact.
\end{itemize}
\end{proposition}
\begin{proof}
(i)$\Leftrightarrow$(ii): See \cite[\S2]{MumfordShimura}  
or \cite[Prop. 17.3.5]{LangeBirkAVII}.

(ii)$\Leftrightarrow$(iii):
On $\Hg(X)(\mathbb{R})$ we have
the conjugation by $h(i)$
\begin{align*}
\Ad h(i):\Hg(X)(\mathbb{R})&\longrightarrow\Hg(X)(\mathbb{R})\\
M&\longmapsto h(i)Mh(i)^{-1},
\end{align*}
which is a Cartan involution (see \cite[\S 2]{DeligneK3}).
If $\Hg(X)$ is commutative, then
$\Ad h(i)$ is just the identity map,
and hence $\Hg(X)(\mathbb{R})$ is compact.

Conversely, if $\Hg(X)(\mathbb{R})$ is compact, then
the identity map is a Cartan involution.
By \cite[Chap.1 Cor. 4.3]{SatakeSymDomain},
any two Cartan involutions of a connected reductive real algebraic group
$\mathcal G$ are conjugate to each other by an inner automorphism of
$\mathcal G$. Hence we have $\Ad h(i)=\Identity$ on $\Hg(X)(\mathbb{R})$.
There are at least two different arguments to finish the proof from the above.

\smallskip
(a) This means that $h(i)$ lies in the centralizer 
$C(\Hg(X)(\mathbb{R}))$ 
of $\Hg(X)(\mathbb{R})$ in $SL(\Gamma_{\mathbb{R}})$.
Suppose we have already shown that one in fact has
\begin{equation}
  \label{eq:ToShow}
  h(i)\in C(\Hg(X))(\mathbb{R}),
\end{equation}
then the whole image $h(S^1)$ lies in $C(\Hg(X))(\mathbb{R})$.
Indeed, since the action $h$ on $\Gamma_{\mathbb{R}}$ is linear, 
we have for any $M\in\Hg(X)$:
\[
h(x+yi)M=h(x)M+h(yi)M=Mh(x)+Mh(yi)=Mh(x+yi).
\]
Since $C(\Hg(X))$ is defined over $\mathbb{Q}$, for $\Hg(X)$ is so,
this shows that $C(\Hg(X))$ contains $\Hg(X)$.
Hence $\Hg(X)$ is commutative.

It remains to prove (\ref{eq:ToShow}). 
Denote the stabilizer of $h(i)$ in $\Hg(X)(\mathbb{C})$ by
\[
H:=\{
M\in \Hg(X)(\mathbb{C})
\mid h(i)M=Mh(i)
\}
\]
Since $H$ is defined by algebraic equations, it is a closed subgroup of
$\Hg(X)(\mathbb{C})$.
Moreover, $H$ contains $\Hg(X)(\mathbb{R})$.
So, from the fact that the $\mathbb{R}$-points 
$\mathcal G(\mathbb{R})$ of any linear algebraic group $\mathcal G$
are Zariski-dense in $\mathcal{G}(\mathbb{C})$,
it follows that
\[
H\supset\overline{\Hg(X)(\mathbb{R})}=\Hg(X)(\mathbb{C}).
\]
Hence 
$H=\Hg(X)(\mathbb{C})$, and $h(i)\in C(\Hg(X)(\mathbb{C}))$.
Since $\mathbb{C}$ is algebraically closed,
we conclude that $C(\Hg(X)(\mathbb{C}))=C(\Hg(X))(\mathbb{C})$,
which implies (\ref{eq:ToShow}).

\smallskip
(b) From
$\Ad h(i)=\Identity$ on $\Hg(X)(\mathbb{R})$
it follows that
for any $\mathbb{Q}$-point $N\in \Hg(X)(\mathbb{Q})$
we have $NI=IN$,
and hence $\Hg(X)(\mathbb{Q})\subset\End_{\mathbb{Q}}(X)$.
On the other hand, we have
\[
\End_{\mathbb{Q}}(X)
=\End(\Gamma_{\mathbb{Q}})^{\Hg(X)}
\]
(see \cite[Prop. 17.3.4]{LangeBirkAVII}),
thus $\Hg(X)(\mathbb{Q})$ is fixed under the conjugation with itself,
whence $\Hg(X)(\mathbb{Q})$ is commutative.

In order to conclude that $\Hg(X)$ is commutative,
we have to show that $\Hg(X)(\mathbb{C})$ is commutative.
By the above, $\Hg(X)(\mathbb{Q})$ lies in the center of 
$\Hg(X)(\mathbb{C})$. 
Recall that the center of a linear algebraic group is a closed subgroup
(see in \cite[Cor. \S 8.2]{HumAlgGr}).
By \cite[Cor. 13.3.9]{SpringerAlgGr} or
\cite[Thm.\,in \S34.4]{HumAlgGr},
$\Hg(X)(\mathbb{Q})$ is Zariski-dense in 
$\Hg(X)(\mathbb{C})$, 
and therefore
\[
\Cent(\Hg(X)(\mathbb{C}))\supset
\overline{\Hg(X)(\mathbb{Q})}
=\Hg(X)(\mathbb{C}),
\]
which implies the commutativity of $\Hg(X)(\mathbb{C})$. 
\end{proof}

Now we prove
\begin{theorem}\label{Prop:RatKMetIffCM}
An abelian variety $X$ is of CM-type
if and only if $X$ admits a constant rational K\"ahler metric.
\end{theorem}
\begin{proof}
$\Leftarrow$:
First suppose $G$ is an arbitrary constant K\"ahler metric on $X$. 
Then for all $z=x+yi\in S^1$ we have
\[
G(h(z)v,h(z)w)=G((x+yI)v,(x+yI)w)
=(x^2+y^2)G(v,w)
=G(v,w),
\]
in other words, $h(S^1)\subset O(G,\mathbb{R})$.
If $G$ is moreover rational, then $O(G)$ is 
an algebraic group defined over $\mathbb{Q}$,
whose $\mathbb{R}$-points
contain $h(S^1)$.
Hence $\Hg(X)$ is an algebraic subgroup of $O(G)$,
and in particular
$\Hg(X)(\mathbb{R})\subset O(G,\mathbb{R})$.
Therefore $\Hg(X)(\mathbb{R})$ is compact,
and $X$ is of CM-type by Proposition \ref{Prop:HgCpcComm}.

$\Rightarrow$:
  If $X$ is of CM-type, then $X$ is isogenous to a product of 
simple abelian varieties of CM-type by Definition \ref{Def:CM}.
So we may assume $X$ is simple.
As mentioned before, $X$ is then isogenous to 
a simple abelian variety of CM-type with a Riemann form $E$ as
in (\ref{Eq:RiemForm}).
It allows us to define the following bilinear 
form on the tangent space:
\[
  G(z,w):=E(z,\beta w)=
\Tr_{K/\mathbb{Q}}(-\beta^2 z\bar w).
\]
We see that $G$ is compatible with $I$ (as $E$ is),
rational, symmetric, and positive definite 
(as $-\beta^2$ is totally positive).
Since the existence of a rational K\"ahler metric 
is preserved under isogeny, this completes the proof.
\end{proof}

\medskip
We shall give an alternative proof of the ``$\Leftarrow$'' direction
of Theorem \ref{Prop:RatKMetIffCM} for a \textit{simple}
abelian variety $X$.
It has the advantage of exhibiting
more clearly how a rational metric
endows $\End_{\mathbb{Q}}(X)$ with additional structures
which force $\End_{\mathbb{Q}}(X)$ to be very ``big''.
Let us first make a reminder of some general facts 
about $\End_{\mathbb{Q}}(X)$.

Due to the presence of the Rosati involution, 
the endomorphism algebra of a simple abelian variety 
must be a division algebra
of finite rank over $\mathbb{Q}$ endowed with a positive anti-involution.
Recall that an involution $f\mapsto f^{\sigma}$ on a division algebra 
$A$ with center $K$ is called positive,
if the quadratic form
\begin{equation}\label{Eq:PosInvo}
\tr_{A|\mathbb{Q}}f^{\sigma}f:=\Tr_{K|\mathbb{Q}}(\tr_{A|K}f^{\sigma}f)
\end{equation}
is positive definite, where $\tr_{A|K}$ denotes the reduced trace of $A$
over $K$,
and $\Tr_{K|\mathbb{Q}}$ denotes the usual trace for the field extension 
$K|\mathbb{Q}$.
Albert gave the classification of such division algebras $A$
(see \cite[Thm 5.5.3, Lemma 5.5.4 and Prop. 5.5.5]{LangeBirkAVII}):
\begin{itemize}
\item[I.] $A=\,$totally real number field, 
left invariant by the positive anti-involution.
\item[II.]$A=\,$totally indefinite quaternion algebra,
there is an element $a\in A$ whose square $a^2$ is in its center $K$
and is totally negative
(i.e.\,is a negative real number under any
embedding $K\hookrightarrow \mathbb{R}$), such that the positive anti-involution
$f\mapsto f^{\sigma}$ is given by $f^{\sigma}=a(\tr_{A|K}f-f) a^{-1}$.
\item[III.]$A=\,$totally definite quaternion algebra, 
and $f\mapsto f^{\sigma}$ is given by $f^{\sigma}=\tr_{A|K}f-f$.
\end{itemize}
The first three algebras are of \textit{the first kind},
i.e. the center is a totally real number field
and coincides with the subfield fixed by the involution.
\begin{itemize}
\item[IV.]The center $K$ of $A$ is a CM-field.
The positive anti-involution restricted to $K$ is the complex conjugation.
\end{itemize}

Let us put $F:=\End_{\mathbb{Q}}(X)$
and denote by $\omega_0$ a (rational) polarization of $X$,
which always exists, since $X$ is algebraic.
The index 0 is to distinguish it from the K\"ahler form $\omega=GI$,
which in general is not rational.
Further we denote by $f\mapsto f'$ the Rosati involution with respect to 
$\omega_0$ and by $G_0$ the K\"ahler metric associated to $\omega_0$.

The presence of a rational K\"ahler metric induces two new structures on $F$:
\begin{itemize}
\item A linear map $\eta\in \End(\Gamma_{\mathbb{Q}})$ determined by
\begin{equation}\label{Eq:EtaDef}
G(\cdot,\cdot)=\omega_0(\eta\cdot,\cdot).
\end{equation}
\item An involution $f\mapsto f^G$ on $\End(\Gamma_{\mathbb{Q}})$
defined by
\[
G(fv,w)=G(v,f^Gw).
\]
\end{itemize}
They have the following properties:

\begin{lemma}\label{Prop:Eta}
  \begin{itemize}
  \item[(i)] $\eta\in F$.
    \item[(ii)] $\eta'=-\eta$.
\item[(iii)] $f\mapsto f^G$ defines a positive anti-involution on $F$.
\item[(iv)] The involution $f\mapsto f^G$ and the Rosati involution are conjugate to 
each other by $\eta$, i.e. $f^G=\eta^{-1} f'\eta,\,\forall f\in F$.
  \end{itemize}
\end{lemma}
\begin{proof}
(i) Since $\omega_0$ and $G$ are compatible with $I$, 
we have for all $v,w\in\Gamma_{\mathbb{Q}}$:
\[
\omega_0(\eta Iv,w)=G(Iv,w)
=-G(v,Iw)
=-\omega_0(\eta v, Iw)
=\omega_0(I\eta v,w),
\]
hence $\eta I=I\eta $, i.e. $\eta \in F$.

(ii) Since 
$\omega_0(\eta v,w)=\omega_0(v,\eta'w)$,
it suffices to show
$\omega_0(\eta v,w)=-\omega_0(v,\eta w)$ for all $v,w\in \Gamma_{\mathbb{Q}}$.
This follows from
\begin{align*}
\omega_0(\eta^{-1} v,w)&=-\omega_0(w,\eta^{-1} v)\\
&=-G(\eta^{-1} w,\eta^{-1} v)\\
&=-G(\eta^{-1} v,\eta^{-1} w)\\
&=-\omega_0(v,\eta^{-1} w).
\end{align*}

(iii) If $f\in F$, i.e. $fI=If$, then $(fI)^G=(If)^G$ and hence
$I^Gf^G=f^GI^G$.
As $I^G=-I$ we get $If^G=f^GI$, i.e. $f^G\in F$.
Next we show that $f\mapsto f^G$
defines a positive anti-involution, i.e.
$\tr_{F|\mathbb{Q}}f^{G}f>0$ for all $f\neq 0\in F$.
Since $F$ acts on $\Gamma_{\mathbb{Q}}$,
one has $\Gamma_{\mathbb{Q}}\cong F^m$,
and $f$ acts on $F^m$ by left multiplication on each component.
On the other hand, the action of $F$ on itself by the left multiplication
has trace $d\cdot\tr_{F|K}f$ over its center $K$, 
where $\tr_{F|K}$ is the reduced trace
and $d^2$ is the degree of $F$ over $K$.
Denote by $\Tr f$ the trace of $f\in F$, when considered
as an endomorphism of $\Gamma_{\mathbb{Q}}$.
Then we have 
\begin{equation}\label{Eq:GPosInvo}
\Tr f=m\cdot d\cdot \Tr_{K|\mathbb{Q}}(\tr_{F|K}f)
=m\cdot\tr_{F|\mathbb{Q}}f,
\quad\forall f\in F.
\end{equation}
In an orthonormal basis with respect to $G$,
$f^G$ is just the transposed matrix of $f$,
hence $\Tr f^Gf>0$, $\forall f\neq 0\in F$.
Then (\ref{Eq:GPosInvo}) implies that the involution 
induced by $G$ is positive.

(iv) This follows immediately from
\[
\omega_0(f'\eta v,w)
=\omega_0(\eta v,fw)
=G(v,fw)
=G(f^Gv,w)
=\omega_0(\eta f^Gv,w)
\]
for all $v,w\in \Gamma_{\mathbb{Q}}$,
which yields $\eta f^G=f'\eta$.
\end{proof}

We shall use Lemma \ref{Prop:Eta} to eliminate the algebras of the first kind,
and then show that the CM-field $K$ of Type IV is
necessarily of maximal rank, i.e. $2g$ over $\mathbb{Q}$,
which implies that $X$ is of CM-type.
\begin{proof}[Second proof of ``$\Leftarrow$'' of Theorem 
\ref{Prop:RatKMetIffCM} for a simple abelian variety $X$]
Type I is already made impossible by (i) and (ii) of Lemma \ref{Prop:Eta}.
On Type III algebras there is a unique
positive anti-involution (see \cite[Prop. 3]{ShimuraPolAV}),
hence $f'=f^G$. 
Then Lemma \ref{Prop:Eta} (iv) implies that $\eta$ lies in $K$,
in contradiction with Lemma \ref{Prop:Eta} (ii).

On Type II algebras, positive anti-involutions are not unique
and they are all of the form given in Albert's classification
mentioned above.
Let us write $f^{\rho}:=\tr_{F|K}f-f$, then
\[
f'=a_1f^{\rho}a_1^{-1}
\quad\mbox{and}\quad
f^G=a_2f^{\rho}a_2^{-1},
\]
for some $a_1$ and $a_2$ in $F$.
Then
$f^G=a_2a_1^{-1}f'a_1a_2^{-1}$,
and, hence, $\eta=\epsilon a_1a_2^{-1}$ for some $\epsilon$ 
in $K$.
On the one hand, $\eta'=-\eta$ by Lemma \ref{Prop:Eta} (ii),
and on the other hand, in view of $a_i^{\rho}=-a_i$
(since $a_i^2\in K$ and $a_i\notin K$), 
we have 
\begin{align*}
  \eta'&=\epsilon a_2^{-1'}a_1'\\
&=\epsilon a_1(a_2^{-1})^{\rho}a_1^{-1}a_1a_1^{\rho}a_1^{-1}\\
&=\epsilon a_1(-a_2^{-1})(-a_1)a_1^{-1}\\
&=\epsilon a_1a_2^{-1}\\
&=\eta.
\end{align*}
A contradiction.

It remains to deal with Type IV algebras.
Recall that the center $K$ of $F$ is a CM-field.
Let us denote by $2m$ its rank over $\mathbb{Q}$
and put $n:=\frac{g}{m}$.
We shall show $n=1$,
which implies that $X$ is of CM-type.
Under $I$ there is the splitting
\begin{equation}\label{Eq:splitting}
\Gamma_{\mathbb{C}}
\overset{\thicksim}{\longrightarrow}
\Gamma_{\mathbb{C}}^{1,0}
\oplus  
\Gamma_{\mathbb{C}}^{0,1}.
\end{equation}
We extend the action of $F$
on $\Gamma_{\mathbb{Q}}$
$\mathbb{R}$-linearly on 
$\Gamma_{\mathbb{R}}$
and denote by $\rho$ its action on $\Gamma_{\mathbb{C}}^{1,0}$
under the isomorphism
$\Gamma_{\mathbb{C}}^{1,0}\cong \Gamma_{\mathbb{R}}$.
Since $K$ is commutative and there is an isomorphism
$\Gamma_{\mathbb{Q}}\cong K^n$,
the action of $K$ on $\Gamma_{\mathbb{Q}}$
diagonalizes on 
$\Gamma_{\mathbb{C}}$,
and the diagonal entries are exactly
$n$ copies of the complete set of $2m$ embeddings of $K$
into $\mathbb{C}$.
The splitting (\ref{Eq:splitting}) 
then implies that $\rho(K)$ even diagonalizes on
$\Gamma_{\mathbb{C}}^{1,0}$,
i.e. there is a complex basis
$\{e_1,\ldots,e_g\}$
of $\Gamma_{\mathbb{R}}$,
with respect to which,
for all $x\in K$ we have $\rho(x)e_l=\rho_l(x)e_l$,
where $\Psi:=\{\rho_1,\ldots,\rho_g\}$
are embeddings of $K$ into $\mathbb{C}$.

We show that $\rho_l$ and $\bar{\rho}_l$
can not both belong to $\Psi$.
Recall that $\eta\in K$ is the element with
$G(\cdot,\cdot)=\omega_0(\eta\cdot,\cdot)$
and $\eta'=-\eta$.
Thus $\rho(\eta)$ has only purely imaginary diagonal entries,
as the Rosati involution restricted to $K$
is just the complex conjugation.
Suppose $\rho_2=\bar{\rho}_1$.
Since $G$ is positive definite on $\Gamma_{\mathbb{R}}$,
we have
\begin{align*}
  0<G(e_1,e_1)
&=\omega_0(\rho(\eta)e_1,e_1)\\
&=\omega_0(\rho_1(\eta)e_1,e_1)\\
&=\Imgn(\rho_1(\eta))\omega_0(Ie_1,e_1)\\
&=\Imgn(\rho_1(\eta))G_0(e_1,e_1),
\end{align*}
and hence $\Imgn(\rho_1(\eta))>0$.
On the other hand,
\[
0<G(e_2,e_2)=\omega_0(\rho_2(\eta)e_2,e_2)
=-\Imgn\rho_1(\eta)G_0(e_2,e_2).
\]
A contradiction.
It follows that $\Psi$ consists exactly of $n$ copies of 
a CM-type $\Phi:=\{\sigma_1,\ldots,\sigma_m\}$
of $K$.

Now we show that this implies 
\begin{equation}\label{eq:IinK}
I\in K\otimes_{\mathbb{Q}}\mathbb{R},
\end{equation}
or equivalently,
$\rho(I)\in\rho(K)\otimes_{\mathbb{Q}}\mathbb{R}$,
where $K\otimes_{\mathbb{Q}}\mathbb{R}$
is considered as a subspace of 
$F\otimes_{\mathbb{Q}}\mathbb{R}$.
Since $\rho(I)$ is just the multiplication by $i$
on $\Gamma_{\mathbb{C}}^{1,0}$
and as we showed that 
$\Psi$ consists of $n$ copies of $\Phi$,
(\ref{eq:IinK})
amounts to show that there is a unique element
$x\in K\otimes_{\mathbb{Q}}\mathbb{R}$,
such that for all $\sigma_l\in\Phi$ we have $\sigma_l(x)=i$.
This is clear due to the isomorphism:
\begin{align*}
K\otimes_{\mathbb{Q}}\mathbb{R}
&\overset{\thicksim}{\longrightarrow}
\mathbb{C}^m\\
a&\longmapsto
(\sigma_1(a),\ldots,\sigma_m(a))^t,
\end{align*}
where we extended $\sigma_l$ $\mathbb{R}$-linearly.

Finally, as $K$ acts by left multiplication on each copy of
$K$ under the isomorphism $\Gamma_{\mathbb{Q}}\cong K^n$
and hence leaves each copy invariant, 
(\ref{eq:IinK})
implies that $I$ leaves each copy of $K\otimes_{\mathbb{Q}}\mathbb{R}$
invariant under the isomorphism
$\Gamma_{\mathbb{R}}\cong (K\otimes_{\mathbb{Q}}\mathbb{R})^n$.
By the simplicity of $X$
we get $n=1$.
This is what we wanted to show.
\end{proof}

\begin{remark}
\textup{From the proof above,
one sees that $\eta$
can be taken as $\beta$ to define the
Riemann form $E$ in (\ref{Eq:RiemForm}).
Moreover, we could also
conclude the proof above by pointing out that 
(\ref{eq:IinK})
means nothing but that the Hodge group
$\Hg(X)$ is contained in $K$,
which implies immediately that 
$\Hg(X)$ is commutative and hence $X$ is of CM-type.}
\end{remark}

%\bigskip
\section{Complex multiplication of the mirror}\label{Sect:CMMirror}
The aim of this section is to show Theorem \ref{Prop:CMMirror}.
It gives
a sufficient condition which ensures that complex multiplication is
transmitted to the mirror partners.
Mirror symmetry for abelian varieties was treated in 
\cite{OrlovMS}.
Some of their results can be rephrased more naturally in terms of 
generalized complex structures,
which were introduced subsequently by 
Hitchin in \cite{HitchinGC}
(see also \cite{GualtieriThesis}, \cite{HuyGC}
and \cite{BBMSGCS} for further works in
generalized complex geometry).
\begin{definition}\label{KapustinGC}
A \textup{generalized complex structure} on a smooth manifold $Y$ is a bundle map
$\mathcal{I}:TY\oplus T^{*}Y\longrightarrow TY\oplus T^{*}Y$
satisfying
\begin{itemize}
\item[(i)] $\mathcal{I}^2=-id$,
\item[(ii)] $\mathcal{I}$ preserves the pseudo-Euclidean metric $q$
on $TY\oplus T^{*}Y$, where
\[
q((a_1,a_2),(b_1,b_2)):=-\langle a_1,b_2\rangle-\langle a_2,b_1\rangle,
\]
\item[(iii)] $\mathcal{I}$ is integrable with respect to the Courant bracket.
\end{itemize}
\end{definition}
\noindent
So the $J_{A\times\hat A}$ and $I_{\omega_A}$ which appeared in 
\cite{OrlovMS} are examples of generalized complex structures,
and the couple
$(J_{A\times\hat A},I_{\omega_A})$
forms a generalized K\"ahler structure,
which is defined as follows 
(see \cite[Chap. 6]{GualtieriThesis} or \cite[\S8]{KGC}).
\begin{definition}\label{Def:gkstr}
  A \textup{generalized K\"ahler structure (GKS)} 
on a smooth manifold $Y$ is a pair $(\mathcal{I},\mathcal{J})$ of
commuting generalized complex structures such that 
$\mathcal G(\cdot,\cdot):=q(\cdot, \mathcal{I}\mathcal{J}\cdot)$ is a positive definite
metric on $TY\oplus T^{\ast}Y$.
\end{definition}
Using GKS, one can define mirror symmetry for a more general class of tori.
\begin{definition}
A \textup{generalized complex torus}
$(\mathbb T,\mathcal{I},\mathcal{J})$
is a real torus $\mathbb T$ endowed with an arbitrary 
GKS $(\mathcal{I},\mathcal{J})$.
\end{definition}
\begin{definition}\label{Def:MS}
Two generalized complex \,\,tori\,\, 
$(\mathbb T=V/\Gamma,\mathcal{I},\mathcal{J})$ \,\,and\,\, 
$(\mathbb T'=V'/\Gamma',\mathcal{I}',\mathcal{J}')$
with arbitrary GKSs
are \textup{mirror} of each other, if there
is a lattice isomorphism
\[
\varphi:\Gamma\oplus\Gamma^*\longrightarrow\Gamma'\oplus\Gamma^{'*},
\]
such that
$q(\cdot,\cdot)=q'(\varphi\cdot,\varphi\cdot)$,
$\mathcal{I}'=\varphi \mathcal{J}\varphi^{-1}$, and
$\mathcal{J}'=\varphi \mathcal{I}\varphi^{-1}.$
We call $\varphi$ a \textup{mirror map}.
We also denote by $\varphi$ its $\mathbb{R}$-linear extension.
\end{definition}
In Theorem \ref{Prop:VAMirror} we will see that 
this notion of mirror symmetry between generalized complex tori
corresponds exactly to the mirror symmetry between the N=2 lattice
vertex algebras associated to them.

Although GKSs provide a more general framework,
a special kind of GKS interests us more particularly.
\begin{definition}\label{Def:IJinduced}
Let $(T,G,B)$ be a complex torus $T$ (in the usual sense) 
with complex structure $I$
and endowed with a constant K\"ahler metric $G$
and a B-field $B\in H^2(T,\mathbb{R})$.
Consider the following pair
\begin{equation}\label{Def:IJ}
\begin{aligned}
  \mathcal{I}:&=\left(
    \begin{matrix}
      1&0\\
      B&1
    \end{matrix}
\right)
\left(
    \begin{matrix}
      I&0\\
      0&-I^{t}
    \end{matrix}
\right)
\left(
    \begin{matrix}
      1&0\\
      -B&1
    \end{matrix}
\right)
=\left(
    \begin{matrix}
      I&0\\
      BI+I^{t}B&-I^t
    \end{matrix}
\right),\\
  \mathcal{J}:&=\left(
    \begin{matrix}
      1&0\\
      B&1
    \end{matrix}
\right)
\left(
    \begin{matrix}
      0&-\omega^{-1}\\
      \omega&0
    \end{matrix}
\right)
\left(
    \begin{matrix}
      1&0\\
      -B&1
    \end{matrix}
\right)
=\left(
    \begin{matrix}
      \omega^{-1}B&-\omega^{-1}\\
      \omega+B\omega^{-1}B&-B\omega^{-1}
    \end{matrix}
\right)\\
&=\left(
    \begin{matrix}
      -IG^{-1}B&IG^{-1}\\
      GI-BIG^{-1}B&BIG^{-1}
    \end{matrix}
\right),
\end{aligned}
\end{equation}
where $\omega$ is the K\"ahler form $G(\cdot,I\cdot)$.
We say that the triple $(T,\mathcal{I},\mathcal{J})$
as defined above
\textup{is induced by} $(T,G,B)$.
\end{definition}
Obviously the above-defined $(T,\mathcal{I},\mathcal{J})$ 
is a generalized complex torus.
Moreover, in this case,
since $\mathcal{I}$ is the B-transform of the usual complex structure
and $\mathcal{J}$ is the B-transform of a symplectic structure,
one can say that a mirror map exchanges
the complex structure and the symplectic structure of the mirror partners
of this kind of generalized complex tori.
Thus, the language of GKS provides a conceptually clean approach 
to mirror symmetry.
We formulate this more precisely in the following 
\begin{definition}
  We say that two complex tori $(T,G,B)$ and $(T',G',B')$
with complex structure $I$ respectively $I'$,
a constant K\"ahler metric $G$ respectively $G'$
and a B-field $B$ respectively $B'$ are 
\textup{mirror partners},
if the generalized complex tori they induce as in Definition 
\ref{Def:IJinduced} are mirror of each other.
\end{definition}

In order to prove Theorem \ref{Prop:CMMirror}
we give the next lemma which first studies the rationality
of the composition $\mathcal{I} \mathcal{J}$ on a
generalized complex torus $(\mathbb T,\mathcal{I},\mathcal{J})$
with an arbitrary GKS,
then links the rationality of $G$ and $B$
to the rationality of $\mathcal{I} \mathcal{J}$
of the GKS they induce.
The rationality of $\mathcal{I} \mathcal{J}$ is defined as follows.
\begin{definition}\label{Def:ratIJ}
Let $(\mathbb T,\mathcal I,\mathcal J)$ be a generalized complex torus
with an arbitrary GKS $(\mathcal{I},\mathcal{J})$.
Let us identify the tangent space 
of $\mathbb{T}$ with its lattice $\Gamma$.
We say that the composition $\mathcal{I} \mathcal{J}$ is 
\textup{defined over} $\mathbb{Q}$,
if $\mathcal{I} \mathcal{J}$ preserves the rational lattices:
\[
\mathcal{I} \mathcal{J}:(\Gamma\oplus\Gamma^*)_{\mathbb{Q}}
\longrightarrow(\Gamma\oplus\Gamma^*)_{\mathbb{Q}}.
\]
\end{definition}

\begin{lemma}\label{Lem:LatticeDecomposition}
Let $(\mathbb T,\mathcal I,\mathcal J)$ be a generalized complex torus
with an arbitrary GKS $(\mathcal{I},\mathcal{J})$
and denote by 
$C_{\pm}\subset(\Gamma\oplus\Gamma^{*})\otimes_{\mathbb{Z}} \mathbb{R}$
the eigenspace of $\mathcal{I}\mathcal{J}$ with eigenvalue $\pm 1$.
Then
\begin{itemize}
\item[(i)] We have the decomposition
\begin{equation}\label{Eq:DecompositionCpm}
(\Gamma\oplus\Gamma^{*})\otimes_{\mathbb{Z}} \mathbb{R}
=C_+\oplus C_-,
\end{equation}
and $C_{\pm}$ is the graph of $\pm id+\mathcal{I} \mathcal{J}$
on $\Gamma_{\mathbb{R}}$.
\item[(ii)] The decomposition (\ref{Eq:DecompositionCpm})
is orthogonal with respect to $q$.
Moreover, $q$ is positive respectively 
negative definite on $C_+$ respectively $C_-$.
\item[(iii)] The decomposition (\ref{Eq:DecompositionCpm})
is defined over $\mathbb{Q}$ if and only if
$\mathcal{I}\mathcal{J}$ is defined over $\mathbb{Q}$.
\item[(iv)] If $(T,\mathcal I,\mathcal J)$ is induced by
$(T,G,B)$, then $C_{\pm}$ is the graph of $\mp G+B$
on $\Gamma_{\mathbb{R}}$,
and (\ref{Eq:DecompositionCpm}) is defined over $\mathbb{Q}$ if and only if
both $G$ and $B$ are rational.
In particular, in this case, $\mathcal{I} \mathcal{J}$ 
is defined over $\mathbb{Q}$ if and only if 
$G$ and $B$ are both rational.
\end{itemize}
\end{lemma}
\begin{proof}
To prove (i), use $\mathcal{I} \mathcal{J}$=$\mathcal{J} \mathcal{I}$.
(ii) follows from $\mathcal{I},\mathcal{J}\in O(q)$
and the requirement that $\mathcal{G} =q \mathcal{I} \mathcal{J}$
is positive definite.
(iii) is a consequence of (i).
To prove (iv) it suffices to note that 
$$  \mathcal{I} \mathcal{J}
\left(
  \begin{matrix}
  1\\
  \mp G+B
  \end{matrix}
\right)
=\pm\left(
  \begin{matrix}
  1\\
  \mp G+B
  \end{matrix}
\right).$$
\end{proof}

The following lemma shows how the lattice of mirror partners
is related to each other.
\begin{lemma}\label{Lem:LatticeDecompositionMirror}
Let $(\mathbb T,\mathcal I,\mathcal J)$ and 
$(\mathbb T',\mathcal I',\mathcal J')$
be mirror generalized complex tori
with arbitrary GKSs
and $\varphi$ a mirror map between them. 
Then
\begin{itemize}
\item[(i)] $\varphi$ respects the decomposition
(\ref{Eq:DecompositionCpm}), i.e.
$\varphi:C_{\pm}\rightarrow C_{\pm}'$.
In particular, $C_+\oplus C_-$ is defined over $\mathbb{Q}$
if and only if $C'_+\oplus C'_-$
is defined over $\mathbb{Q}$.
\item[(ii)] $\varphi$ induces homomorphisms
$\psi_{\pm}:\Gamma_{\mathbb{R}} \overset{\sim}{\rightarrow}\Gamma'_{\mathbb{R}}$
defined by
\[
\varphi(a,\pm a+\mathcal{I} \mathcal{J}a)
=(\psi_{\pm}(a),\pm\psi_{\pm}(a)+\mathcal{I}' \mathcal{J}'\psi_{\pm}(a))
\]
\item[(iii)] 
If 
$(T,\mathcal I,\mathcal J)$ and $(T',\mathcal I',\mathcal J')$
are induced by
$(T,G,B)$ respectively $(T',G',B')$,
then
\begin{itemize}
\item[(a)] $\varphi(a,(\mp G+B)a)
=(\psi_{\pm}a,(\mp G'+B')\psi_{\pm}a)$.
\item[(b)]
$G(a,b)=G'(\psi_{\pm} a,\psi_{\pm} b)$
for all $a,b\in \Gamma_{\mathbb{R}}$.
\item[(c)]
$I'=\psi_+\circ I\circ\psi_+^{-1}
=\psi_-\circ I\circ\psi_-^{-1}$.
\end{itemize}
\end{itemize}
\end{lemma}
\begin{proof}
(i) and (ii) are immediate.
(iii)(a) follows from
Lemma \ref{Lem:LatticeDecomposition} (iv).
For (iii)(b) we make
an explicit calculation
\[
q\left((a,(\mp G+B)a),(b,(\mp G+B)b)\right)
=\pm 2G(a,b)
\quad
\forall a,b\in\Gamma_{\mathbb{R}}.
\]
Then use
$q(\cdot,\cdot)=q'(\varphi\cdot,\varphi\cdot)$ and (iii)(a) to get the claim.
For (iii)(c) we verify the equality for $\psi_-$, the case of $\psi_+$ is 
similar.
Recalling the definition of the generalized K\"ahler structure
 from (\ref{Def:IJ})
we have for any $(a',(G'+B')a')\in C_-'$:
\begin{equation}\label{Eq:I'a'}
\mathcal{I}'\left(
\begin{matrix}
a'\\
(G'+B')a'
\end{matrix}
\right)
=
\left(
\begin{matrix}
I'a'\\
*
\end{matrix}
\right),
\end{equation}
we are only interested in the first component.
Using (iii)(a) the left hand side of (\ref{Eq:I'a'}) is
\begin{equation*}
\varphi \mathcal{J}\varphi^{-1}
\left(
\begin{matrix}
a'\\
(G'+B')a'
\end{matrix}
\right)
=
\left(
\begin{matrix}
\psi_-I\psi_-^{-1}a'\\
(G+B)\psi_- I\psi_-^{-1}a'
\end{matrix}
\right)
\end{equation*}
Comparing with the right hand side of (\ref{Eq:I'a'}),
we obtain (iii)(c).
\end{proof}

As mirror symmetry exchanges complex and K\"ahler structures,
two mirror Calabi-Yau manifolds are in general very different
as complex manifolds.
This remains true for abelian varieties,
but surprisingly the following result shows that
it suffices that the composition $\mathcal{I} \mathcal{J}$ 
is defined over $\mathbb{Q}$ for the mirror to be an isogenous complex torus.
\begin{proposition}\label{Prop:MirrorIsog}
  Let $(T,G,B)$ and $(T',G',B')$ be mirror partners.
If $G$ and $B$ are both rational,
then $T$ and $T'$ are isogenous.
\end{proposition}
\begin{proof}
By Lemma \ref{Lem:LatticeDecomposition} (iv),
$G$ and $B$ are both rational
if and only if 
$\mathcal{I} \mathcal{J}$
is defined over $\mathbb{Q}$.
Then Lemma \ref{Lem:LatticeDecompositionMirror} (i)
implies that $G'$ and $B'$ are rational.
By Lemma \ref{Lem:LatticeDecompositionMirror} (ii)
$\psi_{\pm}$ are then defined over $\mathbb{Q}$.
Finally, from (iii)(c) of the same lemma,
it follows that some integral multiple of $\psi_{\pm}$ 
is actually an isogeny between 
$T$ and $T'$.
\end{proof}
This immediately implies the following 
\begin{theorem}\label{Prop:CMMirror}
  Let $(X,G,B)$
and $(X',G',B')$
be mirror abelian varieties.
Suppose $X$ is of CM-type.
If both $G$ and $B$ are rational,
then $X$ and $X'$ are isogenous.
In particular, $X'$ is also of CM-type.
\end{theorem}
In the next section
we will show that the converse of Theorem \ref{Prop:CMMirror}
however does not hold.
This will also have  some consequence in terms of vertex algebras
(see Corollary \ref{Prop:Summary}).

%\bigskip
\section{An example of mirror abelian varieties of CM-type}
\label{Sect:CounterExp}
In this section we show that the converse of Theorem \ref{Prop:CMMirror}
does not hold.
\begin{proposition}\label{Prop:MirrorNotRat}
There are mirror abelian varieties
$(X,G,B)$ and $(X',G',B')$,
such that $X$ and $X'$ are isogenous and of CM-type,
but neither $\mathcal{I} \mathcal{J}$
nor $\mathcal{I}' \mathcal{J}'$
is defined over $\mathbb{Q}$,
where $(\mathcal{I},\mathcal{J})$
and $(\mathcal{I}',\mathcal{J}')$
denote their induced GKS.
\end{proposition}
We shall eventually construct an explicit example of such mirror pairs,
but we need first some preparation.
\begin{lemma}\label{Lem:cpigamma}
Let
$T\cong \mathbb{C}^g/\Pi \mathbb{Z}^{2g}$ and
$T'\cong \mathbb{C}^g/\Pi' \mathbb{Z}^{2g}$
be complex tori
with period matrix $\Pi$ and $\Pi'$ respectively.
Then $T$ and $T'$ are isogenous
if and only if there is a complex matrix 
$C\in GL(g,\mathbb{C})$ and a rational matrix
$\gamma\in GL(2g,\mathbb{Q})$,
such that
\begin{equation}\label{Eq:PiPi}
\Pi'=C\Pi\gamma.
\end{equation}
In particular, if there is a such matrix $\gamma$, then
$I'=\gamma^{-1}I\gamma$
and if $G$ is a K\"ahler metric on $T$
then $G'(\cdot,\cdot):=G(\gamma\cdot,\gamma\cdot)$
is a K\"ahler metric on $T'$.
\end{lemma}
\begin{proof}
It is easy to see that the rational respectively analytic 
representation of any isogeny provides the matrix 
$\gamma$ respectively $C$.
For the converse, recall
that in general, the rational representation of $I$ of a 
complex torus with period matrix $\Pi$ is 
$I=\left(
\begin{smallmatrix}
\Pi\\
\overline{\Pi}
\end{smallmatrix}
\right)^{-1}
\left(
\begin{smallmatrix}
i \mathbf{1} &0\\
0&-i \mathbf{1} 
\end{smallmatrix}
\right)
\left(
\begin{smallmatrix}
\Pi\\
\overline{\Pi}
\end{smallmatrix}
\right)$.
Replacing 
$\left(
\begin{smallmatrix}
\Pi\\
\overline{\Pi}
\end{smallmatrix}
\right)$
by
$\left(
\begin{smallmatrix}
\Pi'\\
\overline{\Pi}'
\end{smallmatrix}
\right)=
\left(
\begin{smallmatrix}
C &0\\
0&\bar C
\end{smallmatrix}
\right)
\left(
\begin{smallmatrix}
\Pi\\
\overline{\Pi}
\end{smallmatrix}
\right)
\gamma$
for $I'$,
we get $I'=\gamma^{-1}I\gamma$.
It follows that some integral multiple of $\gamma$ 
is an isogeny.
The rest of the claim is obvious.
\end{proof}
Next we give a special form (see (\ref{Eq:PiForm}) below) of the period matrix,
which makes the construction of 
an isogenous mirror easier.
Later we will give an abelian variety of CM-type over a cyclotomic
field, whose period matrix can be written in the form (\ref{Eq:PiForm}).
First a lemma which expresses $I$ explicitly.

\begin{lemma}\label{Lem:ExplicitI}
  Let $\Gamma$ be the lattice of a complex torus $T$
generated by $e_1,\ldots,e_{2g}$.
Suppose that the period matrix $\Pi$ of $T$ in the complex basis 
$\{e_1,\ldots,e_g\}$ 
has the form
$\Pi=\left(
\mathbf{1}\quad T_1+T_2i    
\right)$,
where $T_1$ and $T_2$ are real $g\times g$ matrices,
then in the basis $\{e_1,\ldots,e_{2g}\}$ we have
\[
I=
\left(
\begin{matrix}
-T_1T_2^{-1}&-T_1T_2^{-1}T_1-T_2\\
T_2^{-1}&T_2^{-1}T_1
\end{matrix}
\right).
\]
\end{lemma}
\begin{proof}
The proof is a matter of calculating the matrix
$I=\left(
\begin{smallmatrix}
\Pi\\
\overline{\Pi}
\end{smallmatrix}
\right)^{-1}
\left(
\begin{smallmatrix}
i \mathbf{1} &0\\
0&-i \mathbf{1} 
\end{smallmatrix}
\right)
\left(
\begin{smallmatrix}
\Pi\\
\overline{\Pi}
\end{smallmatrix}
\right)$,
where 
\[
\left(
\begin{matrix}
\Pi\\
\overline{\Pi}
\end{matrix}
\right)^{-1}=\frac{i}{2}\left(
\begin{matrix}
T_1T_2^{-1}-i&-T_1T_2^{-1}-i\\
-T_2^{-1}&T_2^{-1}
\end{matrix}
\right).
\]
\end{proof}

\begin{proposition}\label{Prop:PiAi}
If an abelian variety $X$ has a period matrix
of the form
\begin{equation}\label{Eq:PiForm}
\Pi=\left(
\begin{matrix}
\mathbf{1}&Ai
\end{matrix}
\right),
\end{equation}
with a real matrix $A\in GL(g,\mathbb{R})$,
then by choosing a suitable constant K\"ahler 
metric $G$ and by setting $B=0$,
one can find an isogenous mirror abelian variety 
$(X',G',B')$.
\end{proposition}
\begin{proof}
Suppose that $\Pi$ has the form in (\ref{Eq:PiForm}).
Then by Lemma \ref{Lem:ExplicitI} we have
$I=\left(
\begin{smallmatrix}
0&-A\\
A^{-1}&0
\end{smallmatrix}
\right)$.
Let us choose the metric
\begin{equation*}
G=\left(
\begin{matrix}
-\rho&0\\
0&-A^t\rho A
\end{matrix}
\right),
\end{equation*}
where $\rho$ is a symmetric negative definite matrix 
with integral coefficients.
One verifies that $G$ is compatible with $I$, 
i.e. $I^tG I=G$,
so $G$ is a K\"ahler metric.
Setting $B=0$, then by (\ref{Def:IJ}) $\mathcal{I}$ and $\mathcal{J}$ have the form
\[
\mathcal{I}=
\left(
\begin{matrix}
0&-A&0&0\\
A^{-1}&0&0&0\\
0&0&0&-A^{-1t}\\
0&0&A^t&0
\end{matrix}
\right)
,\quad
\mathcal{J}=
\left(
\begin{matrix}
0&0&0&\rho^{-1}A^{-1t}\\
0&0&-A^{-1}\rho^{-1}&0\\
0&\rho A&0&0\\
-A^t\rho&0&0&0
\end{matrix}
\right).
\]
Further, we choose
\[
C=\rho
\quad\mbox{and}\quad
\gamma=
\left(
\begin{matrix}
\rho^{-1}&0\\
0&\mathbf{1} 
\end{matrix}
\right)
\]
to get a new period matrix
\[
\Pi'=C\Pi\gamma=\left(
\begin{matrix}
\mathbf{1}&\rho Ai
\end{matrix}
\right).
\]
Then by Lemma \ref{Lem:cpigamma} the complex torus 
$X':=\mathbb{C}^g/\Pi' \mathbb{Z}^{2g}$ has complex structure
respectively K\"ahler metric
\begin{equation*}
I'=\gamma^{-1}I\gamma
=\left(
\begin{matrix}
0&-\rho A\\
A^{-1}\rho^{-1}&0
\end{matrix}
\right)
\quad\mbox{ resp. }\quad
G'=\gamma^t G\gamma
=\left(
\begin{matrix}
-\rho^{-1}& 0\\
0&-A^t\rho A
\end{matrix}
\right).
\end{equation*}
Setting $B'=0$ we get
\[
\mathcal{I'}=
\left(
\begin{matrix}
0&-\rho A&0&0\\
A^{-1}\rho^{-1}&0&0&0\\
0&0&0&-\rho^{-1}A^{-1t}\\
0&0&A^t\rho&0
\end{matrix}
\right)
,\quad
\mathcal{J'}=
\left(
\begin{matrix}
0&0&0&A^{-1t}\\
0&0&-A^{-1}&0\\
0&A&0&0\\
-A^t&0&0&0
\end{matrix}
\right).
\]
By easy calculations, one verifies that $(X,G,B)$ and 
$(X',G',B')$
as defined above are mirror partners with the mirror map
\[
\varphi=
\left(
\begin{matrix}
0&0&\mathbf{1}&0\\
0&-\mathbf{1}&0&0\\
\mathbf{1}&0&0&0\\
0&0&0&-\mathbf{1} 
\end{matrix}
\right):
\Gamma\oplus\Gamma^{\ast}
\longrightarrow
\Gamma'\oplus\Gamma^{'\ast}.
\]
Hence, to any abelian variety with period matrix
$\Pi=\left(\mathbf{1}\,\,Ai\right)$ 
by choosing a suitable $G$ and B-field
one can find an isogenous mirror abelian variety.
\end{proof}

In order to obtain a mirror pair of CM-type, 
let us consider the cyclotomic field 
$K=\mathbb{Q}(\xi)$ with 
$\xi^5=1,\xi\neq 1$, which is a CM-field.
Denote by 
\[
w:=e^{\frac{2\pi}{5}i}
=\frac{1}{4}(-1+\sqrt{5})
+\frac{i}{2}\sqrt{\frac{1}{2}(5+\sqrt{5})},
\]
one can write the 
four embeddings of $K$ into $\mathbb{C}$ as
\[
\sigma_k :\xi\longmapsto w^k,
\quad
k=1,\ldots,4.
\]
One has $\sigma_1=\bar{\sigma}_4$ and
$\sigma_2=\bar{\sigma}_3$.
Choose the CM-type $\Phi=\{\sigma_1,\sigma_2\}$, and
choose the lattice to be the ring of integers
$$\Gamma=\mathcal O_K=\mathbb{Z}[\xi].$$
The complex torus 
$X:=\mathbb{C}^2/\Phi(\mathcal O_K)$
is then an abelian variety of CM-type over $K$ 
(see Section \ref{Sect:CMRatMet}).
Let us fix the following generators for $\Gamma$:
\[
\Gamma=\mathbb{Z}\cdot 1
\oplus
\mathbb{Z}(\xi+\xi^{-1})
\oplus
\mathbb{Z}(\xi-\xi^{-1})
\oplus
\mathbb{Z}(\xi^2-\xi^{-2}).
\]
Then under $\Phi$ the lattice is
\[
\Phi(\mathcal O_K)
=\left(
\begin{matrix}
1&w+w^{-1}&w-w^{-1}&w^2-w^{-2}\\
1&w^2+w^{-2}&w^2-w^{-2}&w^4-w^{-4}
\end{matrix}
\right)
\mathbb{Z}
=:
\left(
\begin{matrix}
Z&Ai
\end{matrix}
\right)
\mathbb{Z}.
\]
The left two columns form a real matrix $Z$,
while the right two columns form a purely imaginary matrix
which we write as $Ai$ where $A$ is a real matrix.
Choosing the first two generators to be a complex basis of 
$\Phi(\mathcal O_K)\otimes_{\mathbb{Z}}\mathbb{R}$,
we get the period matrix
\[
\Pi=\left(
\begin{matrix}
\mathbf{1}&Z^{-1}Ai
\end{matrix}
\right)
\]
in the form (\ref{Eq:PiForm}) with the real matrix $Z^{-1}A$.
Together with
\[
G=\left(
\begin{matrix}
-\rho&0\\
0&-A^tZ^{-1t}\rho Z^{-1}A
\end{matrix}
\right)
\quad\mbox{and}\quad
B=0,
\]
the abelian variety $(X,G,B)$
possesses an isogenous mirror
$(X',G',B')$
as constructed in Proposition \ref{Prop:PiAi}.

The last step to get $(\mathcal{I},\mathcal{J})$
such that the composition $\mathcal I\mathcal J$ is not defined over
$\mathbb{Q}$ is to choose an appropriate  $\rho$.
As claimed by 
Lemma \ref{Lem:LatticeDecomposition} (iv),
$\mathcal I\mathcal J$ is defined over
$\mathbb{Q}$
if and only if $G$ and $B$ are both rational.
We set $B=0$ and choose
\[
\rho=\left(
\begin{matrix}
-2&0\\
0&-1
\end{matrix}
\right),
\]
then
\[
G=\left(
\begin{matrix}
2&0&0&0\\
0&1&0&0\\
0&0&5+\frac{2}{\sqrt{5}}&2-\frac{1}{\sqrt{5}}\\
0&0&2-\frac{1}{\sqrt{5}}&3-\frac{2}{\sqrt{5}}\\
\end{matrix}
\right),
\]
which is not rational.
Hence $\mathcal I\mathcal J$ is not defined over $\mathbb{Q}$,
although $X$ is of CM-type and has a mirror
$(X',G',B')$
of CM-type over the same field $K$.
This shows Proposition \ref{Prop:MirrorNotRat}.

%\bigskip
\section{Rationality of lattice vertex algebras,
mirror symmetry and complex multiplication}\label{Sect:VA}
We construct in Section \ref{Subsect:Construction} a vertex algebra 
$V(\Lambda,q,\Lambda_{z})$
in the sense of \cite{KO}.
In the special case where $\Lambda_{z}=\Lambda_{\mathbb R}$
(see (\ref{eq:zzbardecomp}) below),
$V(\Lambda,q,\Lambda_{z})$
shall be the classical lattice vertex algebra
described in \cite[\S5.4]{KacVA}.
We will, therefore, call
$V(\Lambda,q,\Lambda_{z})$
a lattice vertex algebra also in the general case.
In Section \ref{Subset:ToroidalLVA} 
we will see how a generalized complex torus 
$(\mathbb T,\mathcal I,\mathcal J)$
gives rise to such a lattice vertex algebra.
In the appendix we show that it coincides with
the toroidal  vertex algebra in \cite{KO}.
In Section \ref{Subsect:Rat}
we discuss the notion of rationality 
(the precise meaning of this will be recalled).
In Section \ref{Subsect:RatMSCM} we rephrase results of the preceding
sections in terms of lattice vertex algebras,
and answer the question (Q)
posed in the introduction.

\subsection{Construction of lattice vertex algebra
$V(\Lambda,q,\Lambda_{z})$}\label{Subsect:Construction}
We begin with an integral lattice
$(\Lambda,q)$, together with a ($z$,$\bar{z}$)-decomposition 
\begin{equation}\label{eq:zzbardecomp}
\Lambda_{\mathbb{R}}=\Lambda_z\oplus\Lambda_{\bar{z}}
\end{equation}
over $\mathbb{R}$,
which is orthogonal with respect to $q$.
This shall suffice to construct an N=1 vertex algebra.
If moreover we endow the vector space $\Lambda_{\mathbb{R}}$ with an 
almost complex structure $\mathcal{I}$, 
i.e. $\mathcal{I}^2=-\Id$,
then we will get an N=2 vertex algebra
(see \cite[\S 3]{KO} for definitions).
We shall write 
$a=a_z-a_{\bar{z}}$ with $a_z\in \Lambda_z$
and $a_{\bar{z}}\in \Lambda_{\bar{z}}$.
Introduce two copies
$\mathfrak{h}_b=\mathfrak{h}_f=\Lambda_{\mathbb{C}}$
of $\Lambda_{\mathbb{C}}$,
which both inherit the decomposition
\[
\mathfrak{h}_b=\mathfrak{h}_{bz}\oplus\mathfrak{h}_{b \bar{z}}
\quad\mbox{and}\quad 
\mathfrak{h}_f=\mathfrak{h}_{fz}\oplus\mathfrak{h}_{f \bar{z}}.
\]
The affinization is the Lie superalgebra
\[
\hat{\mathfrak{h}}
:=\mathbb{C}[t,t^{-1}]\otimes \mathfrak{h}_b\oplus
t^{\frac{1}{2}}\mathbb{C}[t,t^{-1}]\otimes \mathfrak{h}_f
\oplus \mathbb{C}K
\]
with even
$\mathbb{C}[t,t^{-1}]\otimes \mathfrak{h}_b\oplus\mathbb{C}K$,
and odd
$t^{\frac{1}{2}}\mathbb{C}[t,t^{-1}]\otimes \mathfrak{h}_f$.
The supercommutators are described below by (\ref{Eq:Supercommutator}).
Write $\mathbb{C}[\Lambda]$ for the group algebra of $\Lambda$ over $\mathbb{C}$,
and denote by $e^a,a\in\Lambda$ the basis vectors of
$\mathbb{C}[\Lambda]$.
Furthermore, write
$\mathfrak{h}^<_b:=t^{-1}\mathfrak{h}_b[t^{-1}]$,
and $\mathfrak{h}^<_f:=t^{-\frac{1}{2}}\mathfrak{h}_f[t^{-1}]$.
The space of states is the superspace
\[
V:=\Sym \mathfrak{h}^<_b\otimes\Sym \mathfrak{h}_f^<\otimes \mathbb{C}[\Lambda]
\]
with the vacuum vector $\mathbf{1}=|vac\rangle:=1\otimes 1\otimes 1$.
The parity on $V$ is
$p(s\otimes e^{a})=q(a,a)\mod 2+p(s)$.
The representation $\pi$ of $\hat{\mathfrak{h}}$ 
on $V$ is defined as
\[
\pi:=\pi_1\otimes 1+1\otimes\pi_2,
\]
with $\pi_1$ the representation of $\hat{\mathfrak{h}}$
on 
$\Sym\mathfrak{h}^<_b\otimes\Sym\mathfrak{h}^<_f$
determined by
\begin{align*}
  K&\longmapsto \Id,\\
  h&\longmapsto 0,
\quad h\in \mathfrak{h}_b,\\
n<0,\quad
  t^nh&\longmapsto \mbox{ multiplication by }t^nh,
h\in \mathfrak{h}_b \mbox{ or }h\in \mathfrak{h}_f,\\
n>0,\quad  t^nh&\longmapsto(t^{-s}h'\mapsto n\delta_{n,s}q(h,h')),
\quad h,h'\in \mathfrak{h}_{bz},\\
  t^n\bar h&\longmapsto(t^{-s}\bar h'\mapsto -n\delta_{n,s}q(\bar h,\bar h')),
\quad \bar h,\bar h'\in \mathfrak{h}_{b \bar{z}},\\
  t^nh&\longmapsto(t^{-s}h'\mapsto \delta_{n,s}q(h,h')),
\quad h,h'\in \mathfrak{h}_{fz},\\
  t^n\bar h&\longmapsto(t^{-s}\bar h'\mapsto -\delta_{n,s}q(\bar h,\bar h')),
\quad \bar h,\bar h'\in \mathfrak{h}_{f \bar{z}},
\end{align*}
and $\pi_2$ the representation of $\hat{\mathfrak{h}}$ 
on $\mathbb{C}[\Lambda]$
determined by:
\begin{align*}
  K&\longmapsto 0,\\
  t^nh&\longmapsto(e^a\mapsto \delta_{n,0}q(h,a)e^a),
\quad h\in \mathfrak{h}_{bz},\\
  t^n\bar h&\longmapsto(e^a\mapsto- \delta_{n,0}q(\bar h,a)e^a),
\quad \bar h\in \mathfrak{h}_{b \bar{z}},\\
  t^nh&\longmapsto 0,
\quad h\in \mathfrak{h}_f,\forall n.
\end{align*}
If we write $h_n:=\pi(t^nh)$, then for $h\in\mathfrak{h}_{bz}$,
$\bar h\in\mathfrak{h}_{b\bar z}$,
$f\in\mathfrak{h}_{fz}$,
$\bar f\in\mathfrak{h}_{f\bar z}$,
and $m,n\in \mathbb{Z}$, $r,s\in \mathbb{Z}+\frac{1}{2}$
the supercommutators are 
\begin{equation}\label{Eq:Supercommutator}
\begin{aligned}
&[h_{n},h_{m}']=n\delta_{n,-m}q(h,h'),\quad
[\bar h_{n}, \bar h_{m}']=-n\delta_{n,-m}q(\bar h,\bar h'),\\
&\{f_{r}, f_{s}'\}=\delta_{r,-s}q(f,f'),\quad\quad\,\,\,
\{\bar f_{r},\bar f_{s}'\}=-\delta_{r,-s}q(\bar f,\bar f'),
\end{aligned}
\end{equation}
and all other relations are trivial.
The state-field correspondence maps a homogeneous vector
\begin{equation}\label{Eq:homogeneousVector}
v=h^1_{-s_1}\cdots h^n_{-s_n}
\bar{h}^1_{-\bar{s}_1}\cdots\bar{h}^{\bar{n}}_{-\bar{s}_{\bar{n}}}
f^1_{-r_1}\cdots f^q_{-r_q}
\bar{f}^1_{-\bar{r}_1}\cdots\bar{f}^{\bar{q}}_{-\bar{r}_{\bar{q}}}
\otimes e^{a},
%\quad\in V_{a},
\end{equation}
where $s_i,\bar s_{\bar i}$ are positive integers and
$r_i,\bar r_{\bar i}$ are positive half-integers,
to the field
\begin{equation}\label{eq:sfcorr}
\begin{split}
v(z,\bar z):=&
Y(v,z,\bar z)=\sum_{b\in\Lambda}
\epsilon(a,b)e^a\projection_b
z^{q(a_z,b_z)}\bar{z}^{-q(a_{\bar{z}},b_{\bar{z}})}\\
&\times\exp(-\sum_{n<0}\frac{a_{zn}}{nz^{n}}+\sum_{n<0}\frac{a_{\bar{z}n}}{n\bar{z}^{n}})\\
&\times:\prod^n_{l=1}\frac{\partial^{s_l}H^l(z)}{(s_l-1)!}\prod^{\bar{n}}_{\bar{l}=1}\frac{\partial^{\bar{s}_{\bar{l}}}\bar H^{\bar{l}}(\bar{z})}{(\bar{s}_{\bar{l}}-1)!}\prod^q_{t=1}\frac{\partial^{r_t-\frac{1}{2}}F^t(z)}{(r_t-\frac{1}{2})!}\prod^{\bar{q}}_{\bar{t}=1}\frac{\partial^{\bar{r}_{\bar{t}}-\frac{1}{2}}\bar{F}^{\bar{t}}(\bar{z})}{(\bar{r}_{\bar{t}}-\frac{1}{2})!}:\\
&\times\exp(-\sum_{n>0}\frac{a_{zn}}{nz^{n}}+\sum_{n>0}\frac{a_{\bar{z}n}}{n\bar{z}^{n}}),
\end{split}
\end{equation}
where
$\projection_b$ is the projection onto
$\Sym \mathfrak{h}^<_b\otimes\Sym \mathfrak{h}_f^<\otimes e^b$
and
\begin{alignat*}{2}
\partial H(z)&:=\sum_{m\in \mathbb{Z}}h_mz^{-m-1},
&\qquad
\bar{\partial}\bar H(\bar{z})&:=\sum_{m\in \mathbb{Z}}\bar{h}_m\bar{z}^{-m-1},\\
F(z)&:=\sum_{r\in \mathbb{Z}+\frac{1}{2}}f_rz^{-r-\frac{1}{2}},
&\qquad
\bar{F}(\bar{z})&:=\sum_{r\in \mathbb{Z}+\frac{1}{2}}\bar{f}_r\bar{z}^{-r-\frac{1}{2}},
\end{alignat*}
and the factor $\epsilon(a,b)$ satisfies
the equations (5.4.14) in \cite{KacVA}.
We give a few examples:
\begin{itemize}
\item $v=1\otimes e^0=\mathbf{1}$,
$Y(v,z,\bar{z})=\sum_{b\in\Lambda}\projection_b=id,$
\item $v=h_{-s}\otimes e^0$,
$Y(v,z,\bar{z})=\frac{1}{(s-1)!}\partial^sH(z),$
\item $v=f_{-r}\otimes e^0$,
$Y(v,z,\bar{z})=\frac{1}{(r-\frac{1}{2})!}\partial^{r-\frac{1}{2}}F(z),$
\item $v=1\otimes e^a$,
\begin{equation*}
\begin{split}
Y(v,z,\bar{z})=&\sum_{b\in\Lambda}\epsilon(a,b)e^a\projection_b
z^{q(a_z,b_z)}\bar{z}^{-q(a_{\bar{z}},b_{\bar{z}})}\\
&\times\exp(-\sum_{n<0}\frac{a_{zn}}{nz^{n}}+\sum_{n<0}\frac{a_{\bar{z}n}}{n\bar{z}^{n}})
\exp(-\sum_{n>0}\frac{a_{zn}}{nz^{n}}+\sum_{n>0}\frac{a_{\bar{z}n}}{n\bar{z}^{n}}).
\end{split}
\end{equation*}
\end{itemize}

Next we define the maps $T$ and $\bar{T}$. Let $\{E^i\}$
respectively $\{\bar E^i\}$
be a bosonic basis of $\Lambda_{z}$
respectively $\Lambda_{\bar z}$ and $\{\widetilde E^i\}$
respectively $\{\widetilde{\bar E}^i\}$
be the dual basis with respect to $q$,
i.e. $q(E^i,\widetilde E^j)=\delta^{ij}$
and $q(\bar E^i,\widetilde{\bar E}^j)=-\delta^{ij}$. 
The fermionic bases are denoted by
$\{F^i\}$, $\{\bar F^i\}$ 
and $\{\widetilde F^i\}$, $\{\widetilde{\bar F}^i\}$.
Then
\[
T:=\sum_i\biggl(\sum_{n\geq 0}E^i_n\widetilde{E}^i_{-n-1}
+\sum_{r=\frac{1}{2},\frac{3}{2}\ldots}(r+\frac{1}{2})
F^i_r\widetilde{F}^i_{-r-1}\biggl),\quad
\bar{T}\mbox{ is analogous}.
\]
The superconformal structure is:
\begin{equation*}
    L:=\frac{1}{2}\sum_i (E^i_{-1}\widetilde{E}^i_{-1}-
    F^i_{-\frac{1}{2}}\widetilde{F}^i_{-\frac{3}{2}})
    \otimes e^0,
\quad
\bar{L}\mbox{ is analogous}.
\end{equation*}
\begin{remark}\label{Rem:SpecialKac}
\textup{By inspecting the state-field correspondence
(\ref{eq:sfcorr}) one sees that fields may have non-integral powers 
in $z$ and $\bar{z}$
due to the term 
$z^{q(a_z,b_z)}\bar{z}^{-q(a_{\bar{z}},b_{\bar{z}})}$
(recall that the ($z$,$\bar{z}$)-decomposition 
(\ref{eq:zzbardecomp})
is only required to be defined over $\mathbb{R}$).
However, the difference of the exponents of $z$ and $\bar{z}$ is 
always an integer:
\[
q(a_z,b_z)+q(a_{\bar{z}},b_{\bar{z}})
=q(a,b)\,\in\,\mathbb{Z}.
\]
This implies that
in the special case of a trivial decomposition,
i.e. $\Lambda_{\bar{z}}=0$ or in other words $\Lambda_{\mathbb{R}}=\Lambda_z$,
the lattice vertex algebra is nothing but a conformal
lattice vertex algebra
in the sense of \cite[\S 5.4 \S 5.5]{KacVA} (plus a fermionic part).}
\end{remark}

The N=1 structure is
\[
Q:=\frac{i}{2\sqrt{2}}
\sum_iF^i_{-\frac{1}{2}}\widetilde{E}^i_{-1}
\otimes e^0.
\]
Given an almost complex structure $\mathcal{I}$ 
on the vector space $\Lambda_{\mathbb{R}}$, 
the N=2 structure is denoted by
\begin{align*}
Q^{\pm}:&=\frac{i}{4\sqrt{2}}
\sum_i(F^i_{-\frac{1}{2}}\widetilde{E}^i_{-1}
\pm F^i_{-\frac{1}{2}}\mathcal{I}\widetilde{E}^i_{-1})
\otimes e^0,\\
J:=&-\frac{i}{2} 
\sum_i F^i_{-\frac{1}{2}}\mathcal{I}\widetilde{F}^i_{-\frac{1}{2}}
\otimes e^0.
\end{align*}
and the analogous $\bar{z}$-part.

\subsection{Toroidal lattice vertex algebras}\label{Subset:ToroidalLVA}
Now we explain how tori give rise to lattice vertex algebras.
To a real torus $\mathbb{T}$ together 
with a constant metric $G$ and a B-field,
one can associate a N=1 superconformal lattice vertex algebra
$V(\mathbb{T},G,B)$ by setting
\begin{equation}\label{eq:Lambdaq}
  \begin{aligned}
    &\Lambda=\Gamma\oplus\Gamma^{*},
\quad
q((a_1,a_2),(b_1,b_2)):=-\langle a_1,b_2\rangle-\langle a_2,b_1\rangle\\
\quad\mbox{and}\quad 
&\Lambda_z=\Graph_{\Gamma_{\mathbb{R}}}(-G+B),
\quad
\Lambda_{\bar{z}}=\Graph_{\Gamma_{\mathbb{R}}}(G+B)
  \end{aligned}
\end{equation}
and define the factor 
$\epsilon(a,b):=\exp(i\pi \langle a_1,b_2\rangle)$ for
$a=(a_1,a_2)$ and $b=(b_1,b_2)$
in $\Gamma\oplus\Gamma^*$.

To a generalized complex torus 
$(\mathbb T,\mathcal{I},\mathcal{J})$ with an arbitrary GKS 
$(\mathcal{I},\mathcal{J})$
one can associate two N=2 structures.
Set $(\Lambda,q)$ and $\epsilon(a,b)$
as in (\ref{eq:Lambdaq}) and
\[
\Lambda_z=\Graph_{\Gamma_{\mathbb{R}}}(\Id+\mathcal{I} \mathcal{J}),
\quad
\Lambda_{\bar{z}}=\Graph_{\Gamma_{\mathbb{R}}}(-\Id+\mathcal{I} \mathcal{J}).
\]
Now one can choose either $\mathcal{I}$ or $\mathcal{J}$
to define $Q^{\pm}$ and $J$.
As $\mathcal{I}=\mathcal{J}$ on $\Lambda_z$
and $\mathcal{I}=-\mathcal{J}$ on $\Lambda_{\bar{z}}$
the two N=2 structures are related as follows:
\begin{equation}\label{eq:qjij}
Q^{\pm}_{\mathcal{I}}=Q^{\pm}_{\mathcal{J}},
\quad
\bar Q^{\pm}_{\mathcal{I}}=\bar Q^{\mp}_{\mathcal{J}},
\quad
J_{\mathcal{I}}=J_{\mathcal{J}},
\quad
\bar J_{\mathcal{I}}=-\bar J_{\mathcal{J}}.
\end{equation}
For simplicity, we denote by
$V(\mathbb T,\mathcal{I},\mathcal{J})$ either the N=2 
superconformal lattice vertex algebra 
defined by $\mathcal{I}$ or $\mathcal{J}$.
If additionally, $(T,\mathcal{I},\mathcal{J})$
is induced by $(T,G,B)$, 
where $T$ is a complex torus with complex structure $I$,
$G$ is a K\"ahler metric and $B$ is a B-field,
then we also write
$V(T,G,B)$ for the N=2 lattice vertex algebra $V(T,\mathcal{I},\mathcal{J})$.

In the appendix we prove
\begin{proposition}\label{Prop:IsomoKOVA}
The N=2 superconformal lattice vertex algebra $V(T,G,B)$
is isomorphic to the N=2 superconformal vertex algebra constructed
in \cite{KO}.
\end{proposition}

Adopting the viewpoint of lattice vertex algebras
has the advantage of having a basis-free construction.
Moreover, it is more apparent how the lattice
determines the structure of the vertex algebra.
This leads us to the formulation of rationality
for lattice vertex algebras (see Definition \ref{Def:Rational})
and eventually to prove that the rationality is completely 
determined by the size of the so-called chiral sublattice
$\Lambda_{ch}$ of $\Lambda$ (see Proposition \ref{Prop:VARatMaxRk}).

\subsection{Rationality}
\label{Subsect:Rat}
We phrase the rationality of the lattice vertex algebra 
$V(\Lambda,q,\Lambda_{z})$
in terms of its so-called chiral subalgebra
whose fields contain only integral powers 
in $z$ and $\bar{z}$
(see Definition \ref{Def:Rational}).
Roughly speaking, we call $V(\Lambda,q,\Lambda_{z})$
rational if its chiral subalgebra
is so big that 
$V(\Lambda,q,\Lambda_{z})$
breaks into a finite direct sum of 
irreducible representations of
its chiral subalgebra.

As mentioned in Remark \ref{Rem:SpecialKac}, the term
$z^{q(a_z,b_z)}\bar{z}^{-q(a_{\bar{z}},b_{\bar{z}})}$
may give non-integral powers in $z$ and $\bar{z}$.
However,
the subspace $W\subset V$ of those states which only yield integral powers
is easy to describe.
If we write
\begin{align*}
\Lambda_{ch}:=\{\lambda\in\Lambda
\mid q(\lambda_z,\alpha)\in \mathbb{Z},\forall \alpha\in\Lambda\},
\end{align*}
then we have
\[W=\Sym \mathfrak{h}^<_b\otimes\Sym \mathfrak{h}_f^<
\otimes \mathbb{C}[\Lambda_{ch}].
\]
We call $\Lambda_{ch}$
the \textit{chiral sublattice} of $\Lambda$.
In general, $W$ does not have the structure of a lattice vertex algebra,
as $\Lambda_{ch}$ may be of smaller rank than $\Lambda$.
However, if we require
\begin{itemize}
\item[(a)] $(\Lambda,q)$ is unimodular,
i.e.\,$\Lambda$ is its own dual with respect to $q$, and
\item[(b)] $q$ restricted on $\Lambda_{ch}$ is non-degenerate,
\end{itemize}
then the space
\[
V_{ch}:=\Sym \mathfrak{h}^<_{ch,b}\otimes\Sym \mathfrak{h}_{ch,f}^<
\otimes \mathbb{C}[\Lambda_{ch}],
\]
where $\mathfrak{h}_{ch,b}$ and $\mathfrak{h}_{ch,f}$
are copies of $\Lambda_{ch,\mathbb{C}}$
(see Section \ref{Subsect:Construction}),
does bear the structure of a lattice vertex algebra.
Indeed, if $(\Lambda,q)$ is unimodular,
then the condition
$q(\lambda_z,\alpha)\in \mathbb{Z},\forall \alpha\in\Lambda$
is equivalent to $\lambda_z\in \Lambda$.
Hence denoting
\[
\Lambda_{ch,\mathbb{R}}:=\Lambda_{ch}\otimes_{\mathbb{Z}}\mathbb{R},
\quad
\Lambda_{ch,z}:=\Lambda_{ch,\mathbb{R}}\cap \Lambda_{z},
\quad
\Lambda_{ch,\bar{z}}:=\Lambda_{ch,\mathbb{R}}\cap \Lambda_{\bar{z}},
\]
we get a decomposition
\begin{equation}\label{Eq:DecompCh}
\Lambda_{ch,\mathbb{R}}=\Lambda_{ch,z}\oplus \Lambda_{ch,\bar z}
\end{equation}
which is defined over $\mathbb{Z}$,
i.e.\,$\Lambda_{ch}\cap \Lambda_{z}\subset \Lambda$
and
$\Lambda_{ch,z}=(\Lambda_{ch}\cap \Lambda_{z})\otimes_{\mathbb{Z}}\mathbb{R}$.
Similarly for $\Lambda_{ch,\bar{z}}$.
This gives rise to a lattice vertex algebra 
$V(\Lambda_{ch},q,\Lambda_{ch,z})$ called the 
\textit{chiral subalgebra} of $V(\Lambda,q,\Lambda_{z})$.
Its space of states is exactly $V_{ch}$.
For example, in view of Lemma \ref{Lem:LatticeDecomposition} (ii) and (iv),
the integral lattice $(\Lambda,q)$
and the ($z$,$\bar{z}$)-decomposition 
associated to a torus as defined in (\ref{eq:Lambdaq})
satisfy the conditions (a) and (b), and
hence possess a chiral subalgebra.
Further, the decomposition (\ref{Eq:DecompCh}) allows
us to exhibit the following simple structure of the chiral subalgebra:
\begin{proposition}\label{Prop:ChiralProduct}
Under the conditions (a) and (b) above,
the chiral subalgebra $V(\Lambda_{ch},q,\Lambda_{ch,z})$ 
is the tensor product of two commuting lattice vertex algebras 
in the sense of \cite{KacVA}.
In particular, the fields in $V(\Lambda_{ch},q,\Lambda_{ch,z})$ 
only contain integral powers in $z$ and $\bar{z}$.
\end{proposition}
\begin{proof}
We already showed that the unimodularity of $(\Lambda,q)$
implies that $\Lambda_{ch}\cap \Lambda_{z}$
and $\Lambda_{ch}\cap \Lambda_{\bar{z}}$
are sublattices of $\Lambda$.
Since $\Lambda_z$ is orthogonal to $\Lambda_{\bar{z}}$
with respect to $q$,
the condition (b) implies that $q$ is non-degenerate on 
$\Lambda_z$ and $\Lambda_{\bar{z}}$.
Moreover, since
$\Lambda_{ch,z}=(\Lambda_{ch}\cap \Lambda_{z})\otimes_{\mathbb{Z}}\mathbb{R}$
is contained in $\Lambda_z$,
the ($z$,$\bar{z}$)-decomposition of $\Lambda_{\mathbb{R}}$
induces the trivial decomposition on 
$\Lambda_{ch,z}$, i.e. $\Lambda_{ch,z}$ has no $\bar{z}$-part.
Similarly, $\Lambda_{ch,\bar{z}}$ has no $z$-part.
This gives rise to two lattice vertex algebras
$V(\Lambda_{ch}\cap \Lambda_{z},q,\Lambda_{ch,z})$
and
$V(\Lambda_{ch}\cap \Lambda_{\bar{z}},q,\Lambda_{ch,\bar{z}})$
in the sense of \cite{KacVA}.
Since all the supercommutators between the $z$-
and $\bar{z}$-parts are trivial (see (\ref{Eq:Supercommutator})),
we say that they commute with each other.
Further, if we denote their space of states by $V_{ch,z}$
respectively $V_{ch,\bar{z}}$,
then we have the tensor structure
\[
V_{ch}=V_{ch,z}\otimes V_{ch,\bar{z}}.
\]
Now it is clear that all the fields in the chiral subalgebra 
only contain integral powers in $z$ and $\bar{z}$.
This completes the proof.
\end{proof}

Note that if $\Lambda_{ch}$ is of maximal rank in $\Lambda$,
then we have $V_{ch}=W$,
which means that the chiral subalgebra contains exactly all the fields
with solely integral powers in $z$ and $\bar{z}$.
The maximality of its rank also yields the rationality of 
$V(\Lambda,q,\Lambda_z)$
(see Proposition \ref{Prop:VARatMaxRk}).
In order to give the definition of rationality,
we explain now the module structure of a lattice vertex algebra 
$V(\Lambda,q,\Lambda_z)$
over its chiral subalgebra $V(\Lambda_{ch},q,\Lambda_{ch,z})$.

The state-field correspondence on $V(\Lambda,q,\Lambda_z)$
induces the structure of an $\mathbb{R}^2$-fold algebra on itself,
i.e. for every $(m,\bar m)\in \mathbb{R}^2$, we have an even morphism
\[
V\otimes V\longrightarrow V,
\qquad
a\otimes b\longmapsto a_{(m,\bar m)}b,
\]
where $a_{(m,\bar m)}$ is a coefficient of
$Y(a,z,\bar{z})=\sum_{(m,\bar{m})\in \mathbb{R}^2}
a_{(m,\bar{m})}z^{-m-1}\bar{z}^{-\bar{m}-1}
$.
This structure restricts to a $\mathbb{Z}^2$-fold module structure 
on $V(\Lambda,q,\Lambda_z)$
over its chiral subalgebra (see \cite{RosellenOPEShort} for more details).
\begin{definition}\label{Def:Rational}
A lattice vertex algebra $V(\Lambda,q,\Lambda_z)$ is \textup{rational}
if it decomposes into a finite sum of irreducible modules over its 
chiral subalgebra.
A N=2 lattice vertex algebra is \textup{rational},
if its underlying lattice vertex algebra (without the N=2 structure)
is rational.
\end{definition}
It turns out that the rationality is determined 
by the size of the chiral sublattice.
Indeed, there is an obvious decomposition of $V$. 
For any $\alpha\in\Lambda$, denote by 
$\chi_{\alpha}\subset \mathbb{C}[\Lambda]$ 
the subspace spanned by all vectors
$e^{\alpha+\lambda}$ with $\lambda\in \Lambda_{ch}$.
Clearly, it is independent of the choice of the representant
of $[\alpha]\in\Lambda/\Lambda_{ch}$,
i.e.
$\chi_{\alpha}=\chi_{\alpha'}$
for $[\alpha]=[\alpha']\in\Lambda/\Lambda_{ch}$.
Then
\begin{equation}\label{Eq:DecompXi}
V=\bigoplus_{[\alpha]\in\Lambda/\Lambda_{ch}}
\Sym \mathfrak{h}^<_{b}\otimes\Sym \mathfrak{h}_{f}^<\otimes\chi_{\alpha}
=:\bigoplus_{[\alpha]\in\Lambda/\Lambda_{ch}}V_{\alpha},
\end{equation}
and each $V_{\alpha}$ is a $\mathbb{Z}^2$-fold module over
the chiral subalgebra.
We show
\begin{proposition}\label{Prop:VARatMaxRk}
  A lattice vertex algebra
$V(\Lambda,q,\Lambda_{z})$
is rational if and only if 
$\Lambda_{ch}$ is of maximal rank in $\Lambda$.
\end{proposition}
\begin{proof}
If $\Lambda_{ch}$ is of maximal rank,
i.e. $[\Lambda:\Lambda_{ch}]<\infty$, then the decomposition
(\ref{Eq:DecompXi}) is finite. We show that in this case,
each $V_{\alpha}$ is an irreducible module. Indeed,
as $\Lambda_{ch}$ is of maximal rank,
$V_{ch}$ is isomorphic to $V_{\alpha}$ as vector space, and the action 
of $V_{ch}$ on $V_{\alpha}$ is faithful. 
In view of Proposition \ref{Prop:ChiralProduct}
and \cite[Prop.5.4]{KacVA} we only need to check that
if for some $v\in V_{\alpha}$, we have $(E^i_{-1}\otimes 1)_mv=0$ and $(1\otimes e^a)_mv=0,\forall m,\forall i$ and $\forall a\in\Lambda_{ch}$ 
and similarly for the $\bar{z}$-part, then
$v=0$. 
This is clear by inspecting the explicit expressions of the 
corresponding field of these vectors
(see the examples of fields given in Section \ref{Subsect:Construction}).

Conversely, if $V(\Lambda,q,\Lambda_{z})$ is rational, 
then the sum (\ref{Eq:DecompXi})
must be finite, hence
 $\Lambda_{ch}$ is of maximal rank.
\end{proof}
For the lattice vertex algebra associated to tori,
in view of Lemma \ref{Lem:LatticeDecomposition},
Proposition \ref{Prop:VARatMaxRk} has as consequence the following
\begin{theorem}\label{Thm:VARatGBRat}
(i) The lattice vertex algebra $V(\mathbb T,G,B)$ 
associated to a real torus with
a constant metric $G$ and a B-field
is rational if and only if $G$ and $B$ 
are both rational.

(ii) The N=2 lattice vertex algebra $V(\mathbb T,\mathcal{I},\mathcal{J})$ 
associated to a generalized complex torus with an arbitrary GKS
$(\mathcal{I},\mathcal{J})$
is rational if and only if the composition $\mathcal{I} \mathcal{J}$ 
is defined over $\mathbb{Q}$.

(iii) The N=2 lattice vertex algebra $V(T,G,B)$
associated to a complex torus $T$ endowed with
a constant K\"ahler metric $G$ and a B-field
is rational if and only if $G$ and $B$ 
are both rational.
\end{theorem}
\noindent
In the next section we draw some consequences of Theorem \ref{Thm:VARatGBRat}.

\subsection{Complex multiplication, rationality and mirror symmetry}
\label{Subsect:RatMSCM}
In this paragraph we use Theorem \ref{Thm:VARatGBRat}
to rephrase results of the first part
of the paper in terms of lattice vertex algebras.
Let us start out with an analogue of 
\cite[Thm 5.4]{KO},
which shows that mirror symmetry for generalized 
complex tori can be alternatively expressed in terms of their
associated lattice vertex algebras.
First we recall the definition of mirror symmetry for vertex algebras
from \cite{KO}.
\begin{definition}\label{Def:VAMirror}
  Two N=2 lattice vertex algebras
are \textup{mirror partners} if there is an isomorphism
$f:V\rightarrow V'$ of their space of states, such that
\begin{itemize}
\item[(i)] $f(|vac\rangle)=f(|vac'\rangle)$,
\item[(ii)] $fT=T'f,f\bar T=\bar T'f$,
\item[(iii)] for all $u,v\in V$, we have
$Y'(f(u),z,\bar{z})v=f(Y(u,z,\bar{z})v)$,
\end{itemize}
with the additional property:
\begin{align*}
  &f(Q^+)=Q^{+'},\,f(Q^-)=Q^{-'},\,f(J)=J',\\
  &f(\bar Q^+)=\bar Q^{-'},\,f(\bar Q^-)=\bar Q^{+'},\,f(\bar J)=-\bar J'.
\end{align*}
\end{definition}

\begin{theorem}\label{Prop:VAMirror}
Two generalized complex tori
$(\mathbb T,\mathcal I,\mathcal J)$ and $(\mathbb T',\mathcal I',\mathcal J')$
with arbitrary GKSs
are mirror partners if and only if 
the N=2 lattice vertex algebras
$V(\mathbb T,\mathcal I,\mathcal J)$ and
$V(\mathbb T',\mathcal I',\mathcal J')$ 
are mirror partners.
\end{theorem}
\begin{proof}
Let $\varphi$ be a mirror map between the two generalized
complex tori.
Since $\varphi$ preserves $q$ 
and the decomposition (\ref{Eq:DecompositionCpm})
(see Lemma \ref{Lem:LatticeDecompositionMirror} (i)),
it induces an isomorphism $f$ between the representations $\pi$ of
$\hat{\mathfrak{h}}$ on $V$ and $\pi'$ of $\hat{\mathfrak{h}}'$ on $V'$,
hence $f$ satisfies (i)-(iii) of Definition \ref{Def:VAMirror}.
As $\varphi$ sends $\mathcal{I}\mapsto \mathcal{J}'$,
we have in view of (\ref{eq:qjij})
\begin{align*}
Q^{\pm}_{\mathcal{I}}\mapsto
Q^{\pm'}_{\mathcal{J}'}=Q^{\pm'}_{\mathcal{I}'}
\quad
&J_{\mathcal{I}}\mapsto 
J'_{\mathcal{J}'}=J'_{\mathcal{I}'}\\
\bar Q^{\pm}_{\mathcal{I}}\mapsto
\bar Q^{\pm'}_{\mathcal{J}'}=\bar Q^{\mp'}_{\mathcal{I}'}
\quad
&\bar J_{\mathcal{I}}\mapsto 
\bar J'_{\mathcal{J}'}=-\bar J'_{\mathcal{I}'}
\end{align*}
Similarly for $\mathcal{J}\mapsto \mathcal{I}'$,
hence the N=2 vertex algebra mirror morphism
for both N=2 structures.

Conversely, the isomorphism between the spaces of states induces a bijective map
$\varphi:\Lambda\rightarrow\Lambda'$ of the lattices.
The requirements (i)-(iii) of Definition \ref{Def:VAMirror}
force $\varphi$ to be compatible with
$q$ and $q'$.
Finally, $\varphi$ maps $\mathcal{J}\mapsto \mathcal{I}'$ and 
$\mathcal{I}\mapsto \mathcal{J}'$ because of the N=2 structure of the vertex algebras.
This completes the proof.
\end{proof}

Now we draw a direct consequence of Theorem \ref{Thm:VARatGBRat},
which shows that mirror symmetry has the virtue of letting the
rationality of the lattice vertex algebra
to be transmitted:
\begin{corollary}\label{Prop:MirRatRat}
Suppose $(\mathbb T,\mathcal I,\mathcal J)$ 
and $(\mathbb T',\mathcal I',\mathcal J')$
are mirror generalized complex tori with arbitrary GKSs. Then
the N=2 lattice vertex algebra
$V (\mathbb T,\mathcal I,\mathcal J)$ is rational
if and only if $V (\mathbb T',\mathcal I',\mathcal J')$
is rational.
\end{corollary}
\begin{proof}
From Theorem \ref{Thm:VARatGBRat} (ii)
it follows that $\mathcal{I} \mathcal{J}$
is defined over $\mathbb{Q}$.
Lemma \ref{Lem:LatticeDecomposition} (i)
implies that the decomposition (\ref{Eq:DecompositionCpm})
is defined over $\mathbb{Q}$.
Hence the same holds for $(\mathbb T',\mathcal I',\mathcal J')$
due to Lemma \ref{Lem:LatticeDecompositionMirror} (i).
\end{proof}

Further, by combining Theorem \ref{Thm:VARatGBRat} and 
Theorem \ref{Prop:RatKMetIffCM} one gets
\begin{corollary}\label{Prop:CMVAlgRat}
An abelian variety $X$ is of CM-type
if and only if 
$X$ admits a rational N=2 lattice vertex algebra
$V(X,G,B)$.
\end{corollary}
\begin{proof}
  If $X$ is of CM-type, one can choose $B=0$ together with 
the rational K\"ahler metric claimed by Theorem \ref{Prop:RatKMetIffCM}
to define a rational N=2 lattice vertex algebra
$V(X,G,B)$ in view of Theorem \ref{Thm:VARatGBRat} (iii).
Conversely, the rationality of 
$V(X,G,B)$
forces $G$ to be rational again by Theorem \ref{Thm:VARatGBRat} (iii),
and its N=2 structure means that $G$ is K\"ahler.
Again by Theorem \ref{Prop:RatKMetIffCM}, $X$ is of CM-type.
\end{proof}

Finally, we give an answer to our question
(Q)
on the interplay between complex multiplication,
rationality of the lattice vertex algebra and mirror symmetry 
for abelian varieties. It is actually a reformulation of 
Theorem \ref{Prop:CMMirror} and 
Proposition \ref{Prop:MirrorNotRat}.
\begin{corollary}\label{Prop:Summary}
Let $(X,G,B)$
and $(X',G',B')$ be
mirror abelian varieties.
If the N=2 lattice vertex algebra $V(X,G,B)$ 
is rational,
then $X$ and $X'$ are isogenous and both of CM-type.
Conversely, \,however, 
\,there exist mirror abelian varieties
$(X,G,B)$ and
$(X',G',B') $
such that $X$ and $X'$ are isogenous and both of CM-type,
but neither $V(X,G,B)$ 
nor $V(X',G',B')$ is rational.
\end{corollary}

%----------------------- NEW SECTION -----------------------------------------------------------------

%\bigskip
\section{Appendix: An isomorphism to N=2 SCVA in \cite{KO}}
We show Proposition \ref{Prop:IsomoKOVA}. Recall from \cite[\S3]{KO} the
\begin{definition}\label{Def:valgIsomo}
  Two N=2 superconformal vertex algebras (SCVA)
are \textup{isomorphic} if there is an isomorphism 
$f:V\rightarrow V'$ of their space of states, such that
\begin{itemize}
\item[(i)] $f(|vac\rangle)=f(|vac'\rangle)$
\item[(ii)] $fT=T'f,f\bar T=\bar T'f$
\item[(iii)] For all $u,v\in V$, we have
$Y'(f(u),z,\bar{z})f(v)=f(Y(u,z,\bar{z})v)$
\item[(iv)] $fL=L'$, $fQ^{\pm}=Q^{'\pm}$, $fJ=J'$,
similarly for the $\bar z$-part.
\end{itemize}
\end{definition}
\noindent
Let $(T,G,B)$ be a complex torus.
Recall the charge lattice isomorphism from
\cite{HuybTalkModuli}:
\[
\phi:=\frac{1}{\sqrt{2}}\left(
  \begin{matrix}
    -G-B&1\\
    G-B&1
  \end{matrix}
\right) 
:\Gamma_{\mathbb{R}}\oplus\Gamma_{\mathbb{R}}^{\ast}\longrightarrow 
\Gamma_{\mathbb{R}}^{\ast}\oplus\Gamma_{\mathbb{R}}^{\ast}
\]
with inverse
\begin{equation*}
\phi^{-1}=\frac{1}{\sqrt{2}}\left(
\begin{matrix}
-G^{-1}&G^{-1}\\
1-BG^{-1}&1+BG^{-1}
\end{matrix}
\right).
\end{equation*}
Elementary calculations show that $\phi^{-1}$
is an isometry, i.e.:
\begin{equation}
  \label{eq:qginverse}
  q(\phi^{-1}\cdot,\phi^{-1}\cdot)
=\left(
  \begin{matrix}
    G^{-1}&0\\
    0&-G^{-1}
  \end{matrix}
\right).
\end{equation}
Writing $\phi(a)=:(\phi_1(a),\phi_2(a))$,
the decomposition 
$\Lambda_{\mathbb{R}}
=\Lambda_z\oplus\Lambda_{\bar{z}}
=\Graph_{\Gamma_{\mathbb{R}}}(-G+B)
\oplus\Graph_{\Gamma_{\mathbb{R}}}(G+B)$
from (\ref{Eq:DecompositionCpm})
corresponds to
\[
a\mapsto a_z:=\phi^{-1} \circ\phi_1(a)
\quad\mbox{and}\quad 
a\mapsto a_{\bar{z}}:=-\phi^{-1} \circ\phi_2(a).
\]
Write $f:=\phi^{-1}$. We show that $f$ is the isomorphism
we are looking for.
Here we only give the calculation for the bosonic part.
The fermionic part is similar.

We interpret \cite{KO}'s choice of bases as follows:
$\{\alpha^i\}_{i=1\ldots 2g}$
as the basis of the  \textit{first} component of 
$\Gamma_{\mathbb{R}}^{\ast}\oplus\Gamma_{\mathbb{R}}^{\ast}$
and 
$\{\bar{\alpha}^i\}_{i=1\ldots 2g}$
as the basis of the  \textit{second} component of 
$\Gamma_{\mathbb{R}}^{\ast}\oplus\Gamma_{\mathbb{R}}^{\ast}$.
Set
\[
E^i:=f(\alpha^i)
\quad\mbox{and}\quad 
\bar E^i:=f(\bar{\alpha}^i)
\]
Then for the dual bases, i.e. 
$\{\widetilde{\alpha}^j\}\in\Gamma_{\mathbb{R}}^*\oplus 0$
and
$\{\widetilde{\bar{\alpha}}^j\}\in 0\oplus\Gamma_{\mathbb{R}}^*$
with
$G^{-1}(\alpha^i,\widetilde{\alpha}^j)
=G^{-1}(\bar{\alpha}^i,\widetilde{\bar{\alpha}}^j)
=\delta^{ij}$
set
\[
\widetilde E^i:=f(\widetilde{\alpha}^i)
\quad\mbox{and}\quad 
\widetilde{\bar E}^i:=f(\widetilde{\bar{\alpha}}^i).
\]
Then in view of (\ref{eq:qginverse}) we have
$q(E^i,\widetilde E^j)
=\delta^{ij}$ and
$q(\bar E^i,\widetilde{\bar E}^j)
=-\delta^{ij}$.
On the representations, $f$ induces the correspondence:
\[
\alpha^i_s\mapsto E^i_s,\,\,
\widetilde{\alpha}^i_s\mapsto \widetilde{E}^i_s
\,\,\mbox{for }s\in \mathbb{Z}^*
\quad\mbox{and}\quad
(G^{-1})^{kj}P_k\mapsto E^j_0.
\]
Then the commutators are
\begin{align*}
  [E^i_s,E^j_p]&\overset{(\ref{Eq:Supercommutator})}{=}s\delta_{s,-p}q(E^i,E^j)\\
&\overset{(\ref{eq:qginverse})}{=}s\delta_{s,-p}(G^{-1})^{ij}\\
&=[\alpha^i_s,\alpha^j_p],
\end{align*}
where $s,p\in \mathbb{Z}^*$.
Similarly for the $\bar{z}$-part.

At this stage it is clear that $f$ possesses the properties
(i), (ii) and (iv) of Definition \ref{Def:valgIsomo}.
For (iii) we first translate the notations in \cite{KO} into ours.
For $(w,m)\in \Gamma\oplus\Gamma^{\ast}$:
\begin{align*}
  &P_i(w,m)=\phi_{1i}(w,m)\quad\mbox{is the }i
\mbox{-th coordinate of }\phi_1(w,m)\in\Gamma^{\ast},\\
&\bar P_i(w,m)=\phi_{2i}(w,m)\quad\mbox{is the }i
\mbox{-th coordinate of }\phi_2(w,m)\in\Gamma^{\ast},\\
&k=\phi_1(w,m),\\
&\bar k=\phi_2(w,m).
\end{align*}
Then for $a=(w,m),b=(w',m')\in\Gamma\oplus\Gamma^{\ast}$,
we have again by (\ref{eq:qginverse})
$q(a_z,b_z)=G^{-1}(k,k')$
and $q(a_{\bar{z}},b_{\bar{z} })=-G^{-1}(\bar k,\bar k')$,
and 
\[
\partial^sX(z)
=\partial^{s-1}(G^{-1})^{jk}P_k\frac{1}{z}-\partial^sY^j(z)
\longmapsto
\partial^{s-1}\sum_{m\in \mathbb{Z}}E_mz^{-m-1},
\]
\begin{align*}
  k_jY^j(z)_+=k_j\sum_{m<0}\frac{\alpha^j_m}{mz^m}
\longmapsto
\phi_{1j}(a)\sum_{m<0}\frac{E^j_m}{mz^m}
&=\sum_{m<0}\frac{\phi^{-1}(\sum_j\phi_{1j}(a)\alpha^j)_m}{mz^m}\\
&=\sum_{m<0}\frac{a_{zm}}{mz^m}.
\end{align*}
Compare the state-field correspondence (\ref{eq:sfcorr})
with the one given in \cite{KO},
the isomorphism is then obvious.

%\bigskip\bigskip
\bibliographystyle{amsalpha}

\end{document}